\pdfoutput=1
\documentclass[a4paper, 11pt,reqno]{amsart}
\usepackage[utf8]{inputenc}
\usepackage{amsfonts,extarrows,comment}
\usepackage{bm}
\usepackage{amsmath,amsthm,amssymb,amsfonts,mathrsfs,upgreek,enumerate}
\usepackage{stmaryrd}

\usepackage{tikz}
\usetikzlibrary{intersections,calc}

\usepackage{dsfont}
\usepackage{mathtools}
\usepackage{enumerate}
\usepackage[all]{xy}
\usepackage{dirtytalk}

\usepackage{graphicx}
\usepackage{tikz-cd}
\usepackage{todonotes}
\usepackage{mathtools}
\usepackage{xcolor}
\usepackage{subcaption}
\usepackage{makecell}
\usepackage[margin=1in]{geometry}

\usepackage[colorlinks=true,citecolor=magenta,linkcolor=blue]{hyperref}
\usepackage[justification=centering]{caption}

\newif\ifHideFoot
\HideFoottrue   
\HideFootfalse  

\ifHideFoot

\newcommand{\YHY}[1]{}
\newcommand{\WHY}[1]{}
\newcommand{\CYT}[1]{}

\else 

\newcommand{\marg}[1]{\normalsize{{
			\color{red}\footnote{{\color{blue}#1}}}{\marginpar[\vskip
			-.25cm{\color{red}\hfill$\Rightarrow$\tiny\thefootnote}]{\vskip
				-.2cm{\color{red}$\Leftarrow$\tiny\thefootnote}}}}}

\newcommand{\YHY}[1]{\marg{(YHY) #1}}
\newcommand{\WHY}[1]{\marg{(WHY) #1}}
\newcommand{\CYT}[1]{\marg{(CYT) #1}}

\fi

\usepackage{listings}

\lstdefinestyle{Mathematica}{
    language        =   Mathematica, 
    basicstyle      =   \ttfamily,
    numberstyle     =   \ttfamily,
    keywordstyle    =   \color{blue},
    keywordstyle    =   [2] \color{teal},
    stringstyle     =   \color{magenta},
    commentstyle    =   \color{red}\ttfamily,
    breaklines      =   true,   
    columns         =   fixed,  
    basewidth       =   0.5em,
}

\numberwithin{equation}{section}

1
\numberwithin{subsection}{section}

\allowdisplaybreaks[1]

\usepackage{enumitem}
\setlist[enumerate]{label=\rm{(\roman*)},leftmargin=\parindent,itemindent=\parindent,labelsep=5pt}
\newlist{TFAE}{enumerate}{1}
\setlist[TFAE]{label=\rm{(\alph*)},labelindent=2\parindent}

\newlist{enumeratea}{enumerate}{2}
\setlist[enumeratea]{label=\rm{(\emph{\alph*})}}

\newlist{enumerate1}{enumerate}{2}
\setlist[enumerate1]{label=(\emph{\arabic*})}

\newtheorem*{namedtheorem}{\theoremname}
\newcommand{\theoremname}{testing}

\newtheorem{theorem}{Theorem}[section]
\newtheorem{proposition}[theorem]{Proposition}

\newtheorem{corollary}[theorem]{Corollary}
\newtheorem{lemma}[theorem]{Lemma}

\theoremstyle{definition}
\newtheorem{definition}[theorem]{Definition}

\newtheorem{remark}[theorem]{Remark}
\newtheorem*{remark*}{Remark}

\theoremstyle{remark}

\newcommand\cD{\mathcal{D}}

\newcommand\cF{\mathcal{F}}
\newcommand\cG{\mathcal{G}}
\newcommand\cH{\mathcal{H}}

\newcommand\cM{\mathcal{M}}

\newcommand\cO{\mathcal{O}}

\newcommand\CC{\mathbb{C}}
\newcommand\DD{\mathbb{D}}

\newcommand\PP{\mathbb{P}}
\newcommand\QQ{\mathbb{Q}}

\newcommand\ZZ{\mathbb{Z}}

\newcommand\bB{\mathbf{B}}

\newcommand\bH{\mathbf{H}}

\newcommand\bQ{\mathbf{Q}}

\newcommand\rD{\mathrm{D}}

\newcommand\rH{\mathrm{H}}

\newcommand\rL{\mathrm{L}}

\newcommand\rP{\mathrm{P}}

\newcommand\rT{\mathrm{T}}

\DeclareMathOperator{\D}{D}
\DeclareMathOperator{\id}{id}
\DeclareMathOperator{\Pic}{Pic}

\DeclareMathOperator{\GL}{GL}
\DeclareMathOperator{\PGL}{PGL}
\DeclareMathOperator{\SL}{SL}

\DeclareMathOperator{\Aut}{Aut}

\DeclareMathOperator{\Stab}{Stab}
\DeclareMathOperator{\Fix}{Fix}

\DeclareMathOperator{\im}{\mathrm{Im}}
\DeclareMathOperator{\re}{\mathrm{Re}}

\DeclareMathOperator{\mult}{\mathrm{mult}}
\DeclareMathOperator{\diag}{diag}
\DeclareMathOperator{\proj}{Proj}

\newcommand{\mcolon}{\mathbin{:}}
\newcommand{\defeq}{\vcentcolon=}

\newcommand{\fbrace}[3]{#3\{ #1 \, #3|\, #2 #3\}}

\newcommand{\PTC}{\mathbb{P}{CT}}

\newcommand{\q}{/\!\!/}
\newcommand{\uniL}{\mathbb{I}}
\newcommand{\innerp}{\!\cdot \!}

\title{Compactifications of moduli space of (quasi-)trielliptic K3 surfaces}
\author{Yitao Chen}
\address{Fudan University\\
220 Handan Road\\
200433 Shanghai, China}
\email{19302010030@fudan.edu.cn}

\author{Haoyu Wu}
\address{Shanghai Center for Mathematical Sciences\\
Fudan University\\
2005 Songhu Road\\
200438 Shanghai, China}
\email{hywu18@fudan.edu.cn}

\author{Hanyu Yao}
\address{Fudan University\\
220 Handan Road\\
200433 Shanghai, China}
\email{19300180061@fudan.edu.cn}


\begin{document}
\tikzset{%
    add/.style args={#1 and #2}{
        to path={%
 ($(\tikztostart)!-#1!(\tikztotarget)$)--($(\tikztotarget)!-#2!(\tikztostart)$)%
  \tikztonodes},add/.default={.2 and .2}}
}  

\begin{abstract}
 We study the moduli space $\cF_{T_1}$ of quasi-trielliptic K3 surfaces of type I, whose  general member  is a  smooth bidegree $(2,3)$-hypersurface of $\PP^1\times \PP^2$. Such moduli space plays an important role in the study of the Hassett-Keel-Looijenga program of the  moduli space of degree $8$ quasi-polarized K3 surfaces. 

In this paper,  we consider several natural compactifications of $\cF_{T_1}$, such as the GIT compactification and arithmetic compactifications. 
We give a complete analysis of GIT stability of  $(2,3)$-hypersurfaces and provide a concrete description of the boundary of the GIT compactification.
For the  Baily--Borel compactification of the quasi-trielliptic K3 surfaces, we also compute the configurations of the boundary by classifying certain lattice embeddings. As an application, we show that $(\PP^1\times \PP^2,\epsilon S)$ with small $\epsilon$ is K-stable if $S$ is a K3 surface with at worst ADE singularities. This gives a concrete description of the boundary of the K-stability compactification via the identification of the GIT stability and the K-stability. 
We also discuss the connection between the GIT,  Baily--Borel compactification, and Looijenga’s compactifications by studying the projective models of quasi-trielliptic K3 surfaces. 
\end{abstract}
\maketitle

\section{Introduction}

\subsection{Background}
Let $\cF_{2\ell}$ be the moduli space of primitively quasi-polarized K3 surfaces of degree $2\ell$, i.e., a K3 surface $S$ and a nef line bundle $L$ such that $L^2 =2\ell$ and $c_1(L)\in H^2(S,\ZZ)$ is a primitive class. 
There are natural geometric (partial) compactifications of $\cF_{2\ell}$ constructed by geometric invariant theory.
One can construct the moduli of polarized K3's as the GIT quotient of an open subset of the Hilbert scheme or Chow variety parametrizing $(S, L)$ in $|mL|$ for sufficiently large $m$. 
In the case of low degrees, Mukai showed that general members in $\cF_{2\ell}$ are complete intersections in some homogeneous spaces for $2\ell\leq 22$ and he provided natural GIT compactification of moduli spaces of such K3 surfaces (cf. \cite{Mukai1988}).
On the other hand, by the Global Torelli theorem, $\cF_{2\ell}$ is a locally symmetric variety and hence admits natural arithmetic compactifications such as the Baily--Borel compactification $\cF_{2\ell}^*$ \cite{Baily1966}, Mumford's toroidal compactifications \cite{Ash2010} and Looijenga's semitoric compactifications \cite{Loo03b}. In recent years,  there have been also a series of studies of the compactifications of moduli of K3 surfaces via different methods (cf.~\cite{Brunyate2015,Ascher2019a,Alexeev2020,Laza2016,Hulek2020, Alexeev2019, Alexeev2021}).

A natural question is to investigate the connection between various compactifications.  When  $2\ell=2$, the birational map between $\overline{\cF}_{2}^{\rm GIT}$ and $\cF^\ast_2$ is described by Shah and Looijenga \cite{sha80,looijenga1986new}.
For quartic K3 surfaces, this is the so-called Hassett-Keel-Looijenga program,  proposed by Laza and O'Grady \cite{LO19} for $\cF_4$.  
They  conjectured that  the birational map from the natural GIT compactification to the Baily--Borel compactification could be factorized into a series of elementary birational transformations whose center is the proper transformation of  Shimura subvarieties. 
Recently, this was generalized to larger degrees by Greer-Laza-Li-Si-Tian in \cite{Greer}.
When $2\ell=4$, the program  was confirmed by  Ascher-DeVleming-Liu in \cite{AKDKL2022} by using the HKL program on moduli space of hyperelliptic K3 of degree $4$ (see also  \cite{LO19,LO21}).

\subsection{Moduli of trielliptic K3 surfaces}
In this paper, we investigate the compactification of  moduli space $\cF_{T_n}$ of (quasi)-trielliptic K3 surfaces, of which the general member admits an elliptic fibration with a degree three multisection. 
 The Picard group of a trielliptic K3 surface contains a sublattice $T_n$ whose Gram matrix is 
$$
\begin{array}{c|c|c}
    & C & E    \\ \hline
  C & 2n & 3   \\ \hline
  E & 3 & 0    \\
\end{array} 
.$$
Hence, they are the $T_n$-lattice polarized K3 surface in the sense of \cite{Dolgachev1996}. 
Recently, the families of such K3 surfaces have been studied  by Beauville in  \cite{Beauville2021}.  
Indeed, when $n=1$, the general K3 surface in $\cF_{T_n}$ is a smooth  bidgree $(2,3)$-hypersurface of $\PP^1\times \PP^2$ and there is a natural GIT construction of moduli space of these hypersurfaces  \[
\overline{\cM}^{\rm GIT}_2 \defeq \lvert \cO_{\PP^1\times \PP^2}(2,3) \rvert\, {/\!\!/} \SL_2 \times \SL_3 .
\]
The  moduli space $\cF_{T_1}$  admits a natural morphism to $\cF_8$, whose image is an irreducible Noether-Lefschetz divisor.
 Let $\cF_{T_1}^\ast$ be the Baily--Borel compactification of $\cF_{T_1}$. According to the spirit of HKL program (see \cite{Greer, AKDKL2022}), the natural birational map  $$\overline{\cM}_2^{\rm GIT}\dashrightarrow\cF_{T_1}^\ast$$  has a closed relation to the HKL program of $\cF_8$.  For instance, there will be a one-to-one correspondence between the wall-crossing phenomenon (See more details in Subsection \ref{subsec:HKL}). 
As a start, we  give a detailed description of the GIT compactification and  $\cF_{T_n}^\ast$.  

Let $\pi,\pi'$ be the two projections of $(2,3)$-hypersurface to $\PP^1$ and $\PP^2$ respectively.
We call a curve in a $(2,3)$-hypersurface a horizontal (vertical) line  if it is of the form $\PP^1 \times \{\rm pt\}$ ($\{\rm pt\} \times \PP^1$ respectively).
 Our first main result is

\begin{theorem}\label{thm:main} 
Let $\cM_2$ be the domain of  the birational map $\overline{\cM}_2^{\rm GIT}\dashrightarrow \cF_{T_1}.$ 
Then  the boundary $\overline{\cM}^{\rm GIT}_2 \backslash \cM_2$ consists of $11$ irreducible components $(\alpha)-(r2)$,  whose general member $S$ is described as follows:
\begin{enumerate}
\item The \textbf{strictly semistable} components consist of

$(\alpha)$: $S$ is singular along a vertical line $L$.
And $S$ has a corank $3$ singularity $p$ such that $p \notin \pi(L)$ and $\pi^{\prime}(p) \notin \pi'(L)$.

$(\beta)$: $S$ has two isolated $\widetilde{E}_8$-type singularities $p$ and $q$ such that $\pi(p) \neq \pi(q)$ and $\pi^{\prime}(p) \neq \pi^{\prime}(q)$.
The fibers over $p$ and $q$ are both triple lines with different directions.

$(\gamma)$: $S=(\PP^1\times \PP^1) \cup S'$ for some  $S' \in |\cO_{\PP^1\times \PP^2}(2,2)|$ and $S$ is singular along a horizontal line $C$ such that the intersection of $\PP^1\times \PP^1$ and $C$ is empty.

$(\eta )$: $S$ has two isolated $\widetilde{E}_7$-type singularities $p$ and $q$ such that $\pi(p) \neq \pi(q)$ and $\pi^{\prime}(p) \neq \pi^{\prime}(q)$. The fibers over $\pi(p)$ and $\pi(q)$ contain a line $L_1$ and $L_2$ respectively with $L_1 \neq L_2$. And $\mult_{\pi'(p)}( \pi'(F),L_1)\geq 2$,  $\mult_{\pi'(q)}( \pi'(F),L_2)\geq 2$ for any fiber $F$.

$(\delta)$: $S$ is the union of two $(1,0)$-hypersurfaces and one $(0,3)$-hypersurface.

\item The \textbf{stable} components consist of

$(\zeta)$:   $S$ has an isolated $\widetilde{E}_7$-type singularity $p$ and there is a fiber $F$ such  that $ \pi'(p) \nsubseteq \pi'(F)$.

$(\xi)$:  $S$ has an isolated $\widetilde{E}_8$-type singularity $p$ and the fiber containing $p$ is not a triple line.

$(\theta)$:  $S$ is singular along a section of degree $1$.

$(\phi)$:  $S$ is singular along a section of degree $2$.

$(r1)$:  $S$ is a union of a $(1,1)$-hypersurface and a $(1,3)$-hypersurface.

$(r2)$:  $S$ is a union of a $(0,2)$-hypersurface and a $(2,1)$-hypersurface.
\end{enumerate}

The dimension of the stratum  is  given by 

\[
    \begin{array}{|c|c|c|c|c|c|c|c|c|c|c|c|}
    \hline
      \textrm{\rm Strata}   & \alpha & \beta & \gamma & \eta & \delta & \zeta & \xi & \theta & \phi & r1 & r2 \\ \hline
      \text{\rm dimension}   & 4 & 1 & 2 & 3 & 1 & 10& 7& 8& 5& 13 & 2  \\
    \hline
    \end{array}\]

And the complement of the image of $\cM_2$ in $\cF_{T_1}$ is the union of two irreducible Noether-Lefschetz divisors $\bH_u\cup \bH_h$.
\end{theorem}
Here the NL-divisors $\bH_u$ and $\bH_h$ are the locus of $\cF_{T_1}$ where the $\Pic(S)$ contains a primitive sublattice of the form 
$$
\begin{array}{c|c|c|c}
    & C & E & E'   \\ \hline
  C & 2n & 3  & 1  \\ \hline
  E & 3  & 0  & 1   \\ \hline
  E'& 1  & 1  & 0 \\
\end{array} \quad \text{~ and ~} \quad
\begin{array}{c|c|c|c}
   & C & E & E'   \\ \hline
  C & 2n & 3  & 2  \\ \hline
  E & 3  & 0  & 1   \\ \hline
  E'& 2  & 1  & 0 
\end{array} 
$$ respectively.

\subsection{Connection with K-stable pairs} As proved in \cite{Zhou2021},   there is an open subset of $\cF_{T_1}$ which  admits a K-stability theoretic compactification. 
More precisely, by \cite[Theorem 1.1]{Zhou2021}, there exists  some rational number $0<c<1$ such that one can identify the K-stability of the log Fano pair $(\PP^1\times\PP^2,\epsilon S) $ for $\epsilon \in (0,c)$ with the GIT-stability of $S$ under the action of $\Aut(\PP^1\times \PP^2)$ where $S$ is any $(2,3)$-hypersurface. 
In other words, there is an isomorphism \[
\overline{\cM}^K_2 \defeq M^K_{\PP^1\times \PP^2,\epsilon} \xrightarrow[]{\cong} \overline{\cM}^{\rm GIT}_2
\]
where $M^K_{\PP^1\times \PP^2,\epsilon}$ is a good moduli space parametrizing K-polystable pairs of the form $(\PP^1\times \PP^2,\epsilon S)$. 
As a consequence,  we get
\begin{corollary}     \label{cor:main2}  
There exists a rational number $0 < c< 1$ such that 
the log Fano pair $(\PP^1\times\PP^2,\epsilon S)$ is K-stable for $0<\epsilon<c $ where $S\in |-K_{\PP^1\times \PP^2}|$ is a $(2,3)$-hypersurface with at worst simple singularities (i.e., isolated \rm ADE singularities) or of type $(\zeta),(\xi),(\theta),(\phi)$ or $(r_1), (r_2)$.
\end{corollary}

\subsection{Boundary of Baily--Borel compactification}Following \cite{fr84},  we also compute the boundary of  the Baily--Borel compactification of $\cF_{T_n}$. 

\begin{theorem}[Theorem \ref{thm:Baily--Borel}] \label{thm:main2}
Let $\cF_{T_n}^\ast$ be the Baily--Borel compactification of $\cF_{T_n}$. Then the boundary $\partial\cF_{T_n}:=\cF_{T_n}^\ast-\cF_{T_n}$ is given as follows:
\begin{enumerate}
    \item $n=1$, $\partial \cF_{T_1}$ consists of $14$ modular curves and $2$ points. 
    All curves meet at one point. There are $2$ curves that meet at another point.
    
    \item $n=2$, $\partial \cF_{T_2}$ consists of  $11$ modular curves and $2$ points.
    All curves meet at one point. There are $2$ curves that meet at another point.
    \item $n=3$, $\partial \cF_{T_3}$ consists of $10$  modular curves and $2$ points. 
    All curves meet at one point. There are $2$ curves that meet at another point.
\end{enumerate}
\end{theorem}

In a sequel to this paper, we would like to study the birational transformations between $\overline{\cM}_2^{\rm GIT}$ and $\cF_{T_1}^*$.


 \subsection*{Organization of the paper}
We start in Section \ref{section quasi trielliptic} by constructing $\cF_{T_n}$ the moduli space of the trielliptic K3 surface of the given type.
We also compute the Picard group of $\cF_{T_n}$ using the Noether-Lefschetz theory.
We finish this section by finding out the projective models of trielliptic K3 surfaces.
For $n=1$, generic quasi-trielliptic surfaces are $(2,3)$-hypersurfaces of $\PP^1\times \PP^2$ which gives rise to natural GIT compactification we considered below.
The Sections \ref{section GIT stability}-\ref{section boundary} are devoted to the standard GIT stability analysis of bidegree $(2,3)$-hypersurface. 
Section \ref{section GIT stability} consists of the combinatorics of the unstable and not-properly stable surfaces and we characterize their geometry in  Section \ref{section Geometric interpretation}.
 In section \ref{section minimal orbit}, using Luna's criterion,
 we study the strictly semistable locus of the GIT compactification while in Section \ref{section boundary}, the stable locus with non-simple singularities is discussed.
Together this gives the complete deceptions of the boundary of $\cM_2$ in $\overline{\cM}^{\rm GIT}_2$ and proves the main result Theorem \ref{thm:main}.
 Finally, Baily--Borel compactification and  Looijenga’s compactification of $\cF_{T_n}$ are discussed in section \ref{section Arithmetic compactifications}. 
 In particular, we show that the GIT quotient $\overline{\cM}^{\rm GIT}_2$ is not isomorphic to the Looijenga's compactification associated with $\bH_u \cup \bH_h$.
 We also give a general introduction to the Hassett-Keel-Looijenga program for $\cF_8$ and $\cF_{T_1}$ here.

\subsection*{Acknowledgments}
This project was initiated in the ``Individual Study 2021-2022" of the Talented Undergraduate Program at Fudan University. 
The authors are grateful to Zhiyuan Li for introducing this problem to them.  
The authors also benefit a lot from the discussion with Ruxuan Zhang.  The authors are supported by the NKRD Program of China (No. 2020YFA0713200) and the NSFC General Program (No. 12171090). 
\subsection*{Notations \& Conventions}
\begin{itemize}
    
    \item $\overline{\cM}_2^{\rm GIT}$: the GIT quotient of $(2,3)$-hypersurfaces in $\PP^1 \times \PP^2$.
    \item $\cF_{T_n}$: the moduli space of quasi-trielliptic K3 surfaces of type I,II,III (n=1,2,3) respectively.
    \item $\cM_2$: the domain of the birational map $\overline{\cM}_2^{\rm GIT}\dashrightarrow \cF_{T_1}$.
    \item $\cF_{T_n}^*$: the Baily--Borel compactification of $\cF_{T_n}$.
    \item Let $A_n,\, D_n,\, E_n$ be the negative definite root lattices.
    \item Given an integral lattice $\Lambda$, the discriminant group of $\Lambda$ is $A_{\Lambda} \defeq \Lambda^*/\Lambda$.
    \item Denote by $\uniL_{m,n}$ the set of the even unimodular lattices with signature $(m,n)$.
 \end{itemize}
 We work over $\CC$.

\section{Moduli space of (quasi)-trielliptic K3 surfaces}\label{section quasi trielliptic}

\subsection{(Quasi)-Trielliptic K3 surface} A  K3 surface $S$ is {\bf trielliptic} if it admits an elliptic fibration $\pi \colon S\to \PP^1$ together with a multisection $C$ of degree $3$.  We may say $S$ is trielliptic of type $n$ if $C^2=2n$. In this case, the Picard lattice $\Pic(S)$ contains a primitive sublattice $T_n$ given by
\begin{equation}\label{eq:picard}
   \begin{array}{c|c|c}
 & C & E\\ \hline
  C & 2n & 3   \\ \hline
  E & 3 & 0 \\
\end{array} .
\end{equation}
Up to an isometry, we may assume $n\in \{1,2,3\}$.   
Then $S$ admits a $T_n$-polarization in the sense of \cite{Dolgachev1996}.
Hence we may say a K3 surface is \textbf{quasi-trielliptic} of type I (II, III respectively) if it admits a $T_1$-polarization (respectively $T_2, T_3$-polarization). 

\label{def:lattice}

\subsection{Moduli space of (quasi)-trielliptic K3 surface} \label{subsection:moduli lattice polarized K3}
Let $\cF_{\rT_n}$ be the  coarse moduli space of $\rT_n$-polarized K3 surfaces which parametrizes pairs $(S,\phi) $ where $S$ is a smooth K3 surface, $\phi\mcolon T_n \hookrightarrow \Pic(S)\subset \rH^2(S,\ZZ) $ is a primitive embedding and $\phi(T_n)$ contains a quasi-polarization $L$.
Here two $T_n$-polarized K3 surfaces $(S,\phi)$ and $(S', \phi')$ are equivalent if there exists isomorphism $h\mcolon S'\rightarrow S$ such that $\phi' = h^* \circ \phi$ and  $h^*(L)$ is a quasi-polarization. 
Here the middle cohomology $\rH^2(S,\ZZ)$ is a unimodular even lattice of signature $(3,19)$ under the intersection form $\langle,\rangle$ by cup product. 
For every $S$, we choose a identification of $\rH^2(S,\ZZ) \xrightarrow{\cong} \Lambda = U^2\oplus E_8^3$.
According to \cite{Nikulin1979}, there is a unique primitive embedding up to the automorphism of $\Lambda$
\begin {equation}
    \rT_n\hookrightarrow \Lambda.
\end{equation}
We use $\Sigma_n$ to denote the orthogonal complement of $\rT_n$ in $\Lambda$, which is an even lattice of signature  $(2,18)$.   
Let $\Sigma_{n}^\CC=\Sigma_{n}\otimes \CC$. 
The period domain $\DD$ associated to $\Sigma_{n}$ can be  realized as a connected component of $$\DD^{\pm}:=\{v\in \PP(\Sigma_{n}^\CC)| \left<v,v\right>=0,  \left<v, \bar{v}\right> > 0\}.$$
The monodromy group 
$$\Gamma_{n}=\{g\in O^+(\Sigma_{n})|~g~\hbox{acts trivially on}~\Sigma^\vee_{n}/\Sigma_{n} \}$$
naturally acts on $\DD,$ where $O^+(\Sigma_{n})$ is the identity component of $O(\Sigma_{n})$.
According to the Global Torelli theorem of K3 surfaces (See \cite[Remark 3.4]{Dolgachev1996}), there is an isomorphism $$\cF_{\rT_n}\cong \Gamma_{n}\backslash\DD$$ via the period map. Then $\cF_{\rT_n}$ is a locally Hermitian symmetric variety with only quotient singularities, and hence $\QQ$-factorial.

Note that by the definition of lattice polarization, we may assume the generators $C,E$ of $T_n$ are both effective.
Moreover, up to an automorphism of $\Lambda$, we may assume that $C$ is big and nef (see \cite[Chapter 8, Corollary 2.9]{Huybrechts2016}).

\subsection{Picard group of $\cF_{\rT_n}$} \label{subsec:picard gp}
The Noether-Lefschetz (NL) divisors on $\cF_{\rT_n}$ parametrizes the K3 surfaces in  $\cF_{\rT_n}$ containing additional curve classes.  
According to \cite[Theorem 1]{BLMM17}, the Picard group of $\cF_{\rT_n}$ is generated by NL-divisors.  
Let us recall the construction of some  irreducible NL-divisors on $\cF_{\rT_n}$ and general Heegner divisors. 

\begin{definition}

Let $\beta \in A_{\Sigma_n},m\in \bQ_{<0} $. 
The  Heegner divisor is given by
\begin{equation*}
  \bH_{\beta,m}:= \widetilde{O}^+(\Sigma_n) \setminus \bigcup_{\substack{v\in\Sigma_n^\ast,\, v^2=2m  \\ v\in \beta + \Sigma_n
  }} v^\perp .
\end{equation*}
Note that $ \bH_{\beta,m}= \bH_{-\beta,m}$.

\end{definition} 

In general, the Heegner divisor $\bH_{\beta,m}$ can be non-reduced and reducible, but it can be written as the sum of some irreducible NL divisors which we will introduce below.

\begin{definition} \label{def:NL divisor}
We define  $\bH_u$, $\bH_h$ and $\bH_t$ to be the locus of  K3 surfaces $(S,C,E) \in \cF_{\rT_n}$ such that $\Pic(S)$ contains a divisor class $E'$ satisfying   
\begin{itemize}
    \item   $\bH_u$:  ${E'}^2=0$, $C\innerp {E'}=1$ and $E\innerp {E'}=1$.
    \item $\bH_h$: ${E'}^2=0$, $C\innerp {E'}=2$ and $E\innerp {E'}=1$. 
    \item $\bH_t$: ${E'}^2=0$, $C\innerp {E'}=3$ and $E\innerp {E'}=1$.
\end{itemize}
such that the rank $3$ sublattice of Picard group generated by $\left\{ C,E,{E'} \right\}$ is primitive. 
This is the generalization of primitive NL divisors of moduli space of (quasi)-polarized K3 surfaces.
Following the same proof of \cite[Proposition 1.3]{OGRADY1986}, one can show that these primitive NL divisors are irreducible.

\end{definition}

Now we give the computation of $\Pic(\cF_{T_n})$.
let $L$ be an even lattice of signature $(2,k)$ containing two hyperbolic lattices. 
Let \[
 \rho_L \colon \mathrm{Mp}_2 (\ZZ) \rightarrow \GL(\CC[A_{L}])
\]
be the dual Weil representation of $\mathrm{Mp}_2 (\ZZ)$.
Due to \cite{Bruinier2002} and \cite{BLMM17}, there is an isomorphism\begin{equation} \label{eq:picard-modular form}
 \Pic_{\QQ} (\widetilde{O}^+(L)\setminus \cD_L) \cong \mathrm{Acusp}(\frac{k+2}{2},\rho_L)^{\vee}
\end{equation}
sending Heegner divisor $\bH_{\beta,m}$ to the coefficient function $c_{-m,\beta}$ \[
 \sum_{\tau\in A_L} \sum_{d \in \QQ} c_{d,\tau}q^d e_\tau \mapsto c_{-m,\beta}.
\]
Here $\mathrm{Acusp}(\frac{k+2}{2},\rho_L)$
is the space of almost cusp form of weight $\frac{k+2}{2}$ and type $\rho_L$.
So one can read the relation of Heegner divisors from the the relation between modular forms.
The space of vector-valued modular form can be computed by Raum's method using lattice level Jacobi forms.
For details, we refer to \cite[Section 3, 4]{Peterson2015}.

For the dimension of $\Pic_\QQ(\cF_{T_n})$, one can use the dimension formula of Bruinier \cite{Bruinier2002}:
 \begin{equation}
 \begin{aligned}
 \rho(\cF_{T_n}) \defeq \dim_{\QQ}\Pic_\QQ(\cF_{T_n})
 &= \frac{29}{4} - \frac{1}{12}\re G(2,\Sigma_n)  - \alpha_3(n) - \alpha_4(n)\\
 &\phantom{=} - \frac{1}{9\sqrt{3}} \re \left[ \sqrt{-1}\big(G(1,\Sigma_n) + G(-3,\Sigma_n)\big) \right]
 \end{aligned}
 \end{equation}
where $G(m,L)$ is the generalized quadratic Gauss sum $$G(m,L)=\sum\limits_{\gamma \in A_L} e^{2\pi \sqrt{-1}\frac{m \gamma^2}{2}} $$
and $\alpha_3(n) = \sum_{\gamma \in A_{\Sigma_n}/\pm 1 }  \{-\frac{\gamma^2}{2}\}$, $\alpha_4(n) = \lvert  A_{\Sigma_n}/\pm 1 \rvert$.

Since in our case, the discriminant groups are as follows:
\begin{enumerate}
    \item $n =1, \, A_{\Sigma_1}  = \ZZ\langle \xi_1 \rangle \cong \ZZ/9\ZZ, \, \xi_1^2 = -\frac{10}{9}$,
    \item $ n =2, \, A_{\Sigma_2}  = \ZZ\langle \xi_2 \rangle \cong \ZZ/9\ZZ, \, \xi_2^2 = -\frac{8}{9} $,
    \item $n =3, \, A_{\Sigma_3}  = \ZZ\langle \eta_1 \rangle \times \ZZ\langle \eta_2 \rangle \cong \ZZ/3\ZZ \times \ZZ/3\ZZ, \, \eta_1^2 = -\frac{2}{3}, \, \eta_2^2 =-\frac{4}{3} ,\, \eta_1 \innerp \eta_2 =0$,
\end{enumerate}
 One can compute \[
\alpha_3(n) =\frac{6-n}{3} ,\quad 
\alpha_4(n) =
\begin{cases}
2, & n=1,2 \\
3, & n=3\\
\end{cases}
\]
and \[
G(m,\Sigma_n) = \begin{cases}
3\left(\frac{nm}{9}\right)  , & n=1,2,\ m=1,2 \\
3\sqrt{-3} \left( \frac{-n}{3}\right), & n=1,2,\ m=-3  \\
3\gcd(n,|m|), & n=3,\ m=1,2,-3\\

\end{cases} \quad .
\]
where $\left(\frac{a}{b}\right)$ is the Jacobi symbol.
Therefore we have
\begin{proposition}
 The Picard number of $\cF_{T_n}$ is given by
\[
     \begin{array}{|c|c|c|c|}
        \hline           n & 1 & 2 & 3   \\ \hline
     \rho(\cF_{T_n}) & 3 & 4 & 3 \\ \hline
     \end{array}
\]
and the corresponding Hodge relations are given by
\begin{align*}
    27 \lambda_1 &= 20\bH_{2\xi_1,-\frac{2}{9}} -8 \bH_{\xi_1,-\frac{5}{9}} +\bH_{4\xi_1,-\frac{8}{9}} \\
    &= 21\bH_u -8\bH_h+\bH_t, \\
   108 \lambda_2 &= \bH_{\Bar{0},-1} +130 \bH_{4\xi_2,-\frac{1}{9}} +28\bH_{\xi_2,-\frac{4}{9}} +2\bH_{2\xi_2,-\frac{7}{9}}, \\
    102 \lambda_3 &= \bH_{\Bar{0},-1} +54\bH_{\eta_1,-\frac{1}{3}} + 6\bH_{\eta_2,-\frac{2}{3}}.
\end{align*}
 \end{proposition}

\begin{proof}
The relation follows from the identification \eqref{eq:picard-modular form} and the explicit computation of basis of the almost cusp forms $\mathrm{Acusp}(10, \rho_{\Sigma_n})$ which can be easily computed by \texttt{Sage} package \textbf{weilrep}\footnote{created by Brandon Williams, see \url{https://github.com/btw-47/weilrep}.} using computer.
For examples of quasi-polarized K3 surfaces, see \cite[Section 4.4]{Peterson2015}.
In the case of $\cF_{T_1}$, the second equality follows from the identification of the NL divisors and the Heegner divisors
\[
\bH_u = \bH_{2\xi_1,-\frac{2}{9}} ,\quad \bH_h = \bH_{\xi_1,-\frac{5}{9}}, \quad  \bH_u +  \bH_t = \bH_{4\xi_1,-\frac{8}{9}}.
\]
which can be proved by identifying hyperplanes as in the proof \cite[Section 4,4, Lemma 3]{Maulik2013}.
\end{proof}

\subsection{Projective models of triple-elliptic K3 surfaces}
Let us recall Saint-Donat's classical result of projective models of K3 surfaces. 
\begin{proposition}[see \cite{S-Donat74}]\label{S-Donat}
 Let $(S, L)$ be a smooth \rm K3 surface with a primitive quasi-polarization $L$ of degree $2\ell$ and let $\varphi_L$ be the map defined by $|L|$. 
 Then  one of the following holds:
 \begin{enumerate}
     \item \textup{(Generic case)} $\varphi_L$ birationally maps $S$ to a degree $2\ell$ surface in $\PP^{\ell+1}$. In particular, $\varphi_L$ is a closed embedding when $L$ is ample.
      \item \textup{(Hyperelliptic case)} $\varphi_L$ is a generically $2 : 1$ map and $\varphi_L(S)$ is a smooth rational normal scroll of degree $\ell$, or a cone over a rational normal curve of degree $\ell$.
      Moreover, in this case, $\Pic(S)$ contains a curve class $E'$ satisfying ${E'}^2=0$ and $L\innerp {E'}=2 $ for $L^2 \geq 4$.
      \item \textup{(Unigonal case)} $|L|$ has a fixed component $D$, which is a smooth rational curve. In this case, $\Pic(S)$ contains a curve class $E'$ satisfying ${E'}^2=0$ and $L\innerp {E'}=1$.
 \end{enumerate}
\end{proposition}
The projective models of $\rT_n$-polarized  K3 surfaces are given as below.

\begin{proposition}\label{prop:projective model}
 Let $(S, C, E)$ be a $\rT_n$-quasipolarized \rm K3 surface. 
Consider the rational map $\varphi:S\dashrightarrow |\cO_S(E)|\times |\cO_S(C)|$ and $\overline{\varphi(S)}$ is 
\begin{enumerate}
    \item  a bidegree $(2,3)$- hypersurface of $\PP^1\times \PP^2$ if and only if $(S,C,E)\notin \bH_{u}\cup \bH_h  $ when $n=1$. 
    \item  the complete intersection of two hypersurfaces of bidegree $(1,3), (1,1)$ in $\PP^1\times \PP^3$ if $(S,C,E)\notin  \bH_u \cup \bH_h \cup  \bH_t$ when $n=2$.
 
     \item the intersection of three hypersurfaces of bidegree $(0,3), (1,1), (1,1)$ in $\PP^1\times \PP^4$ if $(S,C,E)\notin \bH_h \cup \bH_t$ when $n=3$.
\end{enumerate}
up to some higher codimension locus.
\end{proposition}
\begin{proof}
According to Saint-Donat's result, the map $\varphi_L$ defined by a primitive quasi-polarization $L$ is a closed embedding after contracting all exceptional $(-2)$ curves if  $L$ is base point free and   $(S,L)$ is not hyperelliptic. 
In our case, we always assume that $C$ is big and nef and $E$ is effective as acknowledged in Subsection \ref{subsection:moduli lattice polarized K3}.
Moreover, one can deduce that $C$ is not unigonal iff $(S,C,E)\notin \bH_{u}$ and when $n\geq 2$, $(S,C)$ is not hyperelliptic iff $(S,C,E)\notin \bH_h$ due to the Hodge Index Theorem.
In the following, we denote by $\widetilde{C}, \widetilde{E}$ the movable part of $C$ and $E$ respectively. 
Hence one can identify $\overline{\varphi(S)}$ with $\Tilde{\varphi}(S)$ where $\Tilde{\varphi}$ is defined by $|\widetilde{E}| \times |\widetilde{C}|$.

For $n=1$, we first show that the intersection matrix given by $(\widetilde{C}, \widetilde{E})$ still equals to $\left(\begin{smallmatrix}
2 & 3\\
3 & 0\end{smallmatrix}\right)$ unless $(S,C,E)\in \bH_u\cup \bH_h$ or in some higher codimension locus. 
We claim that in this case $E$ is necessarily nef hence base point free by \cite[Chapter 2, Proposition 3.10]{Huybrechts2016}.
Assume on the contrary, there exists some irreducible $(-2)$ curve $\Delta$, such that the corresponding intersection matrix is
\[
\begin{array}{c|c|c|c}
       & C & E & \Delta   \\ \hline
     C & 2 & 3 & x  \\       \hline
     E & 3 & 0 & y  \\ \hline
\Delta & x & y & -2 
\end{array}
\]
with $x\geq 0$ and $y<0$.
 It is clear that the only possibilities are $(x,y)=(0,-1)$, $(0,-2)$, $(1,-1)$ or $(2,-1)$ according to the Hodge index theorem. 
 One can see $(S,C,E)\in \bH_h $ if $(x,y)=(1,-1)$ and $(S,C,E)\in \bH_u $ if $(x,y)=(2,-1)$. 
If $(x,y) = (0,-2)$, then the Gram matrix given by basis $\{ C,E, C+\Delta \}$ shows $(S,C,E)\in \bH_h$.
If $(x,y)=(0,-1)$, one sees that $(E-\Delta)$ is nef unless $S$ lies in some higher codimension locus.
Thus $(E-\Delta)$ is base point free as $(E-\Delta)^2=0$ and we have $\widetilde{E}= E-\Delta$.
Thus $(S,C,E)\notin \bH_u\cup \bH_h$, the Gram matrix given by $(\widetilde{C}=C, \widetilde{E})$ still equals to $\left(\begin{smallmatrix}
2 & 3\\
3 & 0\end{smallmatrix}\right)$
Set $L= C+\widetilde{E}$.
Moreover, since the sublattice spanned by $\{C,\widetilde{E}\}$ is obtained from $T_1$ by a reflection respect to $\Delta$,  we see $L$ is not hyperelliptic or unigonal by assumption.
One can see that  $s\circ \tilde{\varphi} =  \tilde{\varphi}_{|L|}$ is a closed embedding after contracting all exceptional $(-2)$ curves, where $s$ is the Segre embedding.
Then $\tilde{\varphi} $ is birational to its image  which is a $(2,3)$-hypersurface by the adjunction formula.

Conversely, for $(S,C,E)\in \bH_u$, one find that $\Delta\defeq E- E'$ is effective and  makes $E$ not nef.
And it's easy to find that $\widetilde{E} = E'$ and $\widetilde{C} = 2E'$ with $h^0(S,\cO_S(2E')) = 3$ by \cite[Proposition 2.6 and 2.7.2]{S-Donat74}.
Then the image of the projection $\overline{\varphi(S)} \to |\widetilde{C}|$ is a rational normal curve.
Hence $\overline{\varphi(S)}$ can not be a $(2,3)$-hypersurface.
For $(S,C,E)\in \bH_h$, like the above case we have $\widetilde{E} = E'$ and $\overline{\varphi(S)} $ is equal to the image induced by $|C| \times |E'|$ given by intersection matrix $\left(\begin{smallmatrix}
2 & 2\\
2 & 0\end{smallmatrix}\right)$.
Notice that the projection map $\overline{\varphi(S)} \to |C| \cong \PP^2$ restricted to $\varphi(\widetilde{E})$ has degree $2$.
So the projection $\overline{\varphi(S)} \to \PP^2$ is generically injective, hence $\varphi$ is a generically $2$ to $1$ map and $\overline{\varphi(S)}$ is not a $(2,3)$-hypersurface.

For $n=2$ and $3$, one can show that $E$ is nef unless $(S,C,E) \in \bH_u \cup \bH_h \cup \bH_t $ using the similar analysis as $n=1$ case.
Note that $H_u = \emptyset $ for $n=3$.
Hence when $S$ does not lies in the union of NL divisors, we have $E$ is base point free and $\varphi$ is a morphism.
Note that $\varphi_{|C|}$ is a closed embedding after contracting all exceptional $(-2)$ curves.
Then $\varphi=\varphi_{|C| \times |E|}$ is also a closed embedding. 

Next, we analyse the projective model $\varphi(S)$ for $n=2$ and $3$. 
For $n=2$, note that the composition with the Segre embedding $s\circ \varphi\colon S \to \PP^7$ is induced by the line bundle $\cO_S(L)$ and the dimension of linear system $|C+E|$ is $6$.
Then $s \circ \overline{\varphi(S)}$ is contained in a hyperplane and so $\overline{\varphi(S)}$ is contained in a bidegree $(1,1)$-hypersurface $X_{1,1}$ of $\PP^1\times \PP^3$.  
Using the Lefschetz hyperplane theorem and adjunction formula, one can see that the divisor class of $S$ in $X_{1,1}$ is $(1,3)$. Note that $$H^1(I_{X_{1,1}}(1,3)) = H^1(\cO_{\PP^1 \times \PP^3}(0,2))=0,$$ where $I_{X_{1,1}}$ is the ideal sheaf of $X_{1,1}$.
We get a surjection 
$H^0(\cO_{\PP^1 \times \PP^3}(1,3)) \twoheadrightarrow H^0(\cO_{X_{1,1}}(1,3))$. Thus $\overline{\varphi(S)}$ is the complete intersection of two hypersurfaces of bidegree $(1,3)$, $(1,1)$ in $\PP^1\times \PP^3$. 

For $n=3$, similarly, $s\circ \varphi\colon S \to \PP^1 \times \PP^4 \xrightarrow{s} \PP^9$ is induced by the line bundle $\cO_S(L)$ where $s$ is the Segre embedding and the dimension of linear system $|C+E|$ is $7$, one can see that $s \circ \overline{\varphi(S)}$ is contained in the intersection of two hyperplanes.
So $\overline{\varphi(S)}$ is contained in the intersection of two bidegree $(1,1)$-hypersurfaces $X_1 
$ 
and 
$X_2 
$ of $\PP^1\times \PP^4$.
We denote $X_1 \cap X_2$ by $Y$. One can see that if a $(1,1)$- hypersurface is singular, the only possibility is that it is reducible. Thus $X_1$ and $X_2$ are smooth. Otherwise, $S$ is contained in a $\{\textrm{pt}\} \times \PP^4$ or a $\PP^1 \times \PP^3$, which is impossible. 
By checking the Jacobian of $Y$, one can conclude that $Y$ is also smooth.
Note that the divisor class of $S$ in $Y$ is $(0,3)$ by adjunction. 
We deduce that the divisor class $(0,3)$ of $Y$ comes from a bidegree $(0,3)$-hypersurface of $\PP^1\times \PP^4$ since $H^1(\cO_{X_1}(-Y)(0,3))$ and $H^1(\cO_{\PP^1 \times \PP^4}(-X_1)(0,3))$ vanish. 
Therefore, $\overline{\varphi (S)}$ is the complete intersection of three hypersurfaces of bidegree $(0,3)$, $(1,1)$, $(1,1)$ in $\PP^1\times \PP^4$.


\end{proof}

\section{Stability of bidegree \texorpdfstring{$(2,3)$}{(2,3)}-hypersurfaces in \texorpdfstring{$\PP^1\times \PP^2$}{P1*P2}}\label{section GIT stability}

\subsection{Numerical criteria}

Using Hilbert-Mumford's numerical criteria \cite[Thm. 2.1]{Mumford1994}, we have: A bidegree $(2,3)$-hypersurfaces in $\PP^1\times \PP^2$ is stable (resp. semistable) if and only if $\mu(f,\lambda) > 0$ (resp. $\geq 0$) for all 
one parameter subgroups $\lambda$ of $\SL_{2}\times \SL_{3}$, where
$\mu(f,\lambda)$ is the numerical weight introduced by Hilbert and Mumford.

As is customary, a one parameter subgroup (1-PS) of $\SL_{2}\times \SL_{3}$ can be diagonalized as
$$\lambda \colon t\in \CC^{\ast} \rightarrow \diag(t^{a},t^{-a},t^b,t^c,t^{-b-c})$$ for some $a,b,c \in \ZZ$. We call such $\lambda$ a normalized 1-PS of $\SL_{2}\times \SL_{3}$ if $a\geq 0$ and $b\geq c \geq -b-c$.

Let $\lambda$ be a normalized 1-PS. Then the weight of a  monomial $x_0^u x_1^{2-u} y_0^v y_1^w y_2^{3-v-w}$ with respect to $\lambda$ is
$$au-a(2-u) +bv+cw+ (-b-c)(3-v-w).$$
If we denote by $M^{\circleddash} (\lambda)$ (resp. $M^{-} (\lambda)$) the set of monomials of bidegree $(2,3)$ which have non-positive (resp. negative) weight with respect to $\lambda$, one can easily compute the maximal
subsets $M^{\circleddash} (\lambda)$ (resp. $M^{-} (\lambda)$) listed in the next subsection.

\begin{table}[ht]
    \centering
     \begin{tabular}{|c|c|c|c|}
 \hline
 Cases & $1$-PS  &Maximal monomials  &  Invariant \\
 \hline
 N1 & $\lambda'_1 = (3,-3,2,2,-4)$  & $x_1^2 y_0^3, x_0 x_1 y_0^2 y_2, x_0^2 y_0y_2^2$  &  $(\alpha)$ \\
 \hline
 N2 & $\lambda'_2 = (3,-3,2,0,-2)$  & $x_1^2 y_0^3, x_0 x_1 y_0 y_1 y_2, x_0^2 y_2^3$ &  $(\beta)$\\
 \hline
 N3 & $\lambda'_3 = (3,-3,4,-2,-2)$  & $x_1^2 y_0^2y_1, x_0 x_1 y_0 y_1^2, x_0^2 y_1^3$  & $(\alpha)$
 \\  
 \hline
 N4 & $\lambda'_4 = (0,0,2,-1,-1)$  & $ x_0^2 y_0y_1^2$ & $(\gamma)$
 \\
 \hline
 N5 & $\lambda'_5 = (1,-1,2,0,-2)$  & $x_1^2 y_0^2 y_2, x_0 x_1 y_0 y_1 y_2, x_0^2 y_0 y_2^2 ,x_0^2y_1^2y_2$ & $(\eta)$
 \\
 \hline
 N6 & $\lambda'_6 = (0,0,1,1,-2)$  & $ x_0^2 y_0^2 y_2$ & $(\gamma)$
 \\
 \hline
 N7 & $\lambda'_7 = (1,-1,0,0,0)$  & $ x_0 x_1 y_0^3$ & $(\delta)$
 \\
 \hline
\end{tabular}
    \caption{Not properly stable}
    \label{tab:NPS}
\end{table}

\begin{table}[ht]
    \centering
     \begin{tabular}{|c|c|l|c|}
 \hline
 Cases & $1$-PS  & Description (roughly) & inclusion\\
 \hline
 U1 & $\lambda_1 = (5,-5,3,-1,-2)$  & reducible  & S1 \\
 \hline
 U2 & $\lambda_2=(4,-4,2,1,-3)$  & singular along a vertical line & S2, S3  \\
 \hline
 U3 & $\lambda_3 = (4,-4,4,-1,-3)$  & corank 3 isolated & S2, S3 ($\star$)\\  
 \hline
 U4 & $\lambda_4 = (3,-3,4,-1,-3)$  & corank 3 isolated & S3\\
 \hline
 U5 &  $\lambda_5 = (1,-1,3,-1,-2) $ & singular along a horizontal line & S3 
\\
 \hline
 U6 & $\lambda_6 = (1,-1,5,-1,-4)$  & singular along a horizontal line & S4\\
 \hline
 U7 & $\lambda_7 = (2,-2,3,1,-4)$  & singular along a vertical line & S5 \\
 \hline
\end{tabular}
    \caption{Unstable}
    \label{tab:my_label}
\end{table}

\subsection{Maximal subsets for not properly stable points}\label{numerical}

\begin{lemma} \label{lem:ss orbit}
For any normalized $1$-\rm PS $\lambda$, $M^{\circleddash} (\lambda)$ is contained in one of $M^{\circleddash}(\lambda'_i) $ in table \ref{tab:NPS}.
The surface $S$ is not properly stable if its defining polynomial is one of the following: 
\begin{itemize}
    \item[\rm(N1)] $
     f(x_0,x_1,y_0,y_1,y_2) = x_1^2 c(y_0,y_1,y_2) + x_0x_1 y_2 q(y_0,y_1) + x_0^2 y_2^2 l(y_0, y_1 , y_2). $
     
     \item[\rm(N2)] $f(x_0,x_1,y_0,y_1,y_2) = 
    x_1^2 c_0(y_0,y_1,y_2) +  x_0 x_1 [c_1(y_1,y_2) +  y_0 y_2 l(y_1, y_2)] + \mu x_0^2 y_2^3 $.
    
    \item[\rm(N3)] $
    f(x_0,x_1,y_0,y_1,y_2) = 
    x_1^2 [c_1(y_1,y_2) + y_0 q_1(y_1,y_2) + y_0^2 l(y_1,y_2)] \\
    + x_0 x_1 [c_2(y_1,y_2) + y_0 q_2(y_1,y_2)] + x_0^2 c_0(y_1,y_2)$.
    
    \item[\rm(N4)] $
    f(x_0,x_1,y_0,y_1,y_2) = x_1^2 [c_0(y_1,y_2) + y_0 q_0(y_1,y_2)] + x_0 x_1 [c_1(y_1,y_2) + y_0 q_1(y_1,y_2)] \\ + x_0^2 [c_2(y_1,y_2) + y_0 q_2(y_1,y_2)]$.
    
    \item[\rm(N5)] $f(x_0,x_1,y_0,y_1,y_2) = x_1^2[c_0(y_1,y_2) + y_0 q(y_1,y_2)+ \mu y_0^2 y_2] + x_0x_1[c_1(y_1,y_2) + y_0y_2 l(y_1, y_2)] 
    \\+ x_0^2 y_2 (q(y_1,y_2) + \nu y_0y_2)$.
    
    \item[\rm(N6)] $f(x_0,x_1,y_0,y_1,y_2) = y_2 \cdot q(x_0,x_1,y_0,y_1,y_2)$.
    
    \item[\rm(N7)] $f(x_0,x_1,y_0,y_1,y_2) = x_1 [x_0 c_0(y_0,y_1,y_2) + x_1 c_1(y_0,y_1,y_2)]$.
\end{itemize}
\end{lemma}

 Using the destabilizing 1-PS, the invariant part of equation (N1)-(N7) are given as follows: 
  \begin{align}
(\alpha) :&\; f = x_1^2 c(y_0,y_1)+x_0x_1 q(y_0,y_1) y_2+ x_0^2 l(y_0,y_1) y_2^2.\label{alpha}\\
(\beta) :&\; f= a x_1^2 y_0^3 + x_0x_1 ( by_1^3+cy_0y_1y_2)+ dx_0^2 y_2^3.\label{beta} \\
(\gamma) :&\; f = x_1^2 q_1(y_0,y_1)y_2+x_0x_1 q_2(y_0,y_1)y_2+ x_0^2 q_3(y_0,y_1)y_2. \label{gamma}\\
(\eta) :&\; f= x_1^2 (a_1y_1^2y_0+b_1y_0^2 y_2)+x_0x_1(a_2y_0y_1y_2+b_2y_1^3)+ x_0^2
    (a_3y_0y_2^2+b_3y_1^2y_2). \label{eta}\\
(\delta) :&\;f= x_0 x_1 c(y_0,y_1,y_2). \label{delta}
\end{align}


Similarly, we can get maximal subsets for unstable points.

\begin{lemma}\label{lem:us orbit}
For any normalized $1$-\rm PS $\lambda$, $M^{-} (\lambda)$ is contained in one of $M^- (\lambda_i)$ in table \ref{tab:my_label}. 
The surface $S$ is unstable if its defining polynomial is one of the following: 
\begin{itemize}
    \item[\rm(U1)] $
    f(x_0,x_1,y_0,y_1,y_2) = x_1^2 c(y_0,y_1,y_2) + x_0x_1[c_0(y_1,y_2) + \mu y_0 y_2^2].$
    \item[\rm(U2)] $
    f(x_0,x_1,y_0,y_1,y_2) =x_1^2 c_0(y_0,y_1,y_2) + x_0x_1 (y_2^2 l(y_0,y_1,y_2) + \mu y_2 y_1^2) + \nu x_0^2 y_2^3.$
    \item[\rm(U3)] $
    f(x_0,x_1,y_0,y_1,y_2) = x_1^2 [ c_0(y_1,y_2) + y_0  q(y_1,y_2) + y_0^2 l(y_1,y_2) ] +x_0 x_1 [c_1(y_1,y_2)+ \mu y_0  y_2^2] \\+ \nu x_0^2 y_2^3.$

    \item[\rm(U4)] $f(x_0,x_1,y_0,y_1,y_2) = x_1^2 [c_0(y_1,y_2) + y_0 q(y_1,y_2) + \mu y_0^2 y_2 ]  + x_0x_1 [c_1(y_1,y_2) +\nu y_0 y_2^2]  \\+x_0^2 y_2^2 l(y_1,y_2)$.

    \item[\rm(U5)] $ f(x_0,x_1,y_0,y_1,y_2) = x_1^2 [c_0(y_1,y_2) + y_0 q_1(y_1,y_2))]
    + x_0 x_1 [c_1(y_1,y_2)+ \mu y_0 y_2^2]
    +x_0^2 c_2(y_1,y_2)$.

    \item[\rm(U6)]
    $f(x_0,x_1,y_0,y_1,y_2) = x_1^2 [c_0(y_1,y_2) 
    + y_0y_2 l(y_1, y_2)]
    +x_0 x_1 [c_1(y_1,y_2)+ \mu y_0 y_2^2]  \\+x_0^2 [     c_2(y_1,y_2) + \nu y_0 y_2^2].$
    
    \item[\rm(U7)] $f(x_0,x_1,y_0,y_1,y_2) = x_1^2[c(y_1,y_2)+ y_0 y_2 l_0(y_1,y_2) + \mu y_0^2y_2] 
    + x_0x_1 (y_2^2 l_1(y_0,y_1,y_2)+ \nu y_2y_1^2) \\+ x_0^2 y_2^2 l_2(y_0,y_1,y_2)$.
    
\end{itemize} 
\end{lemma}

\section{Geometric interpretation of stability}\label{section Geometric interpretation}

Let us first give some notations and conventions.
\subsection{Additional Conventions} 
Given a bidegree $(2,3)$-hypersurface  $S \subseteq \PP^1\times \PP^2$, it admits an elliptic fibration 
$$\pi: S\to \PP^1$$
via the first projection and all the fibers we consider are respect to $\pi$. 
On the other hand, the second projection 
$$\pi':S \to\PP^2$$
is a double cover branching along a sextic curve denoted by $\bB(S)$.
In particular,  we call a curve  $C\subseteq S$ 
\begin{itemize}
    \item \emph{a vertical line}, if $C$ has the form $\{\mathrm{pt}\}\times \PP^1$.
    
    \item \emph{a horizontal line} if $C$ has the form $\PP^1 \times \{\text{pt}\}$.
    
    \item \emph{a section of degree $d$} if $\pi(C) = \PP^1$ and $\pi^{-1}(p)\innerp C  =d $ for any $p\in \PP^1$.
\end{itemize}

In the rest of this paper, We use the terminology of the corank of the hypersurface singularities 
as in \cite{Arnold2012} and \cite{Laza2009}. 
\begin{definition}
Let $0\in \CC^n$ be a hypersurface singularity given by an equation $f(z_1,\ldots,z_n)=0$. The corank of $0$ is $n$ minus the rank of the Hessian of $f(z_1,\ldots,z_n)$ at $0$. 
\end{definition}
Moreover, let $f(u, v, w)$ be an analytic function in $\CC [[u,v,w]]$ whose leading term defines an isolated singularity at the origin. 
We are concerned with the following analytic types of isolated hypersurface singularities:
\begin{itemize}
  \item Simple singularities: $A_n (n\geq 1)$, $D_n (n\geq 4)$ and $E_r(r=6,7,8)$.
  \item Simple elliptic singularities $\widetilde{E}_r (r=6,7,8)$:
  \begin{itemize}
      \item $\widetilde{E}_6$: $f= u^3 + v^3 +w^3 + uvw$.
      \item $\widetilde{E}_7$: $f= u^2 + v^4 +w^4 + uvw$.
      \item $\widetilde{E}_8$: $f= u^2 + v^3 +w^6 + uvw$.
  \end{itemize}
  
\end{itemize}
\subsection{Geometry of not properly stable surface}

\begin{theorem}\label{thm:semistable}
   A bidegree $(2,3)$-hypersurface $S$ in $\PP^1\times \PP^2$ is not properly stable if and only if  one of the following conditions holds
   \begin{itemize}
       \item[\rm (N1)] $S$ is singular along a vertical line;
       \item[\rm (N2)] $S$ contains a singularity $p$ of  at least  $\widetilde{E}_8$-type, and the fiber over $\pi(p)$ is a triple line.
       \item[\rm (N3)] $S$ contains at least a singularity of corank $3$; 
       \item[\rm (N4)] $S$ is singular along a horizontal line;
       \item[\rm (N5)] $S$ contains a corank $2$ singularity $p$ of at least $\widetilde{E}_7$-type, such that the fiber over $\pi(p)$ contains a line $L$, and $\mult_{\pi'(p)}(\pi'(F),L)\geq 2$ for any fiber $F$.
       \item[\rm (N6)] $S=(\PP^1\times \PP^1) \cup S'$ for some  $S' \in |\cO_{\PP^1\times \PP^2}(2,2)|$;
       \item[\rm (N7)] $S=\PP^2\cup S'$ for some  $S' \in |\cO_{\PP^1\times \PP^2}(1,3)|$.
   \end{itemize}
\end{theorem}
\begin{proof}
According to Lemma \ref{lem:ss orbit}, it suffices to find the geometric characterizations of type {\rm (N1)-(N7)}. 
We give detailed proof for types (N2) and (N5) here.
The other cases are similar and relatively simple so we omit them.

If $S$ is of type (N2), then the equation of $S$ is given by 
\[
 x_1^2 c_0(y_0,y_1,y_2) +  x_0 x_1 [c_1(y_1,y_2) +  y_0 y_2 l(y_1, y_2)] + \mu x_0^2 y_2^3=0.
\]
   One can assume that the coefficient of $x_1^2y_0^3$ and $x_0^2y_2^3$ is nonzero. 
   Otherwise, the equation will degenerate to type (N3) and (N6). \
   Set $p=(1,0,1,0,0)$, then the fiber over $(1,0)$ is the triple line $3L\colon \{y_2^3=0\}$.
Letting $x_0=y_0=1$, the affine equation near $p$ is 
\[
  x_1^2 +y_2^3 + x_1^2 l(y_1,y_2) + x_1y_2 l(y_1, y_2) + x_1^2q(y_1,y_2) +x_1c_1(y_1,y_2)=0.
\]
It is clear that the weight on variables $x_1$, $y_2$, $y_1$ is $(\frac{1}{2},\frac{1}{3},\frac{1}{6})$. Thus $p$ is at least $\widetilde{E}_8$-type.

Conversely, we take $p= (1,0,1,0,0)$ to be the isolated singularity of at least $\widetilde{E}_8$-type.
Up to coordinate change, one can assume the fiber over $\pi(p)$ is the triple line $3L\colon \{y_2^3=0\}$, 
then the equation of $S$ can be written as
\[
  x_1^2 c_0(y_0,y_1,y_2) + x_0 x_1 [c_1(y_1,y_2) + y_0 q_1(y_1,y_2)] + \mu x_0^2 y_2^3=0.
\]
If the coefficient of $x_0x_1y_0y_1^2$ is nonzero, one can see that the weight of $p$ has weight no worse than $(\frac{1}{2},\frac{1}{3},\frac{1}{4})$, which contradicts the condition that $p$ is at least $\widetilde{E}_8$-type.

If $S$ is of type (N5), then the equation of $S$ is given by 
\[
\begin{split}
   x_1^2[c_0(y_1,y_2) + y_0 q(y_1,y_2)+ \mu y_0^2 y_2] + x_0x_1[c_1(y_1,y_2) + y_0y_2 l(y_1, y_2)] 
    \\ + x_0^2 y_2 (q(y_1,y_2) + \nu y_0y_2)=0.
\end{split}
\]
One can assume that the coefficients of $x_1^2y_0^2y_2$ and $x_0^2y_0y_2^2$ are nonzero. Otherwise the equation will degenerate to type (N3) and (N4).
Set $p=(1,0,1,0,0)$, then the fiber over $(1,0)$ contains the line $\{y_2=0\}$. 
And $\mult_{(1,0,0)}(\pi'(F), y_2)\geq 2$ for any fiber $F$ since there are no $y_0^3$, $y_0^2y_1$ terms.
One can check that $p$ is at least $\widetilde{E}_7$-type.

Conversely, we take $p= (1,0,1,0,0)$ to be the isolated corank $2$ singularity of at least  $\widetilde{E}_7$-type.
Up to coordinate change, one can assume that the fiber over $(1,0)$ contains the line $\{y_2=0\}$, then the equation of $S$ can be written as
\[
  x_1^2 c_0(y_0,y_1,y_2) + x_0 x_1 [c_1(y_1,y_2) + y_0 q_1(y_1,y_2) +  y_0^2 y_2] + x_0^2 y_2 (q(y_1,y_2) + \nu y_0y_2)=0.
\]
One can deduce that the coefficients of $x_1^2y_0^3$ and $x_1^2y_0^2y_1$ are zero since 
$\mult_{\pi'(p)}(\pi'(F),y_2)\geq 2$ for any fiber $F$. 
And the coefficient of $x_0x_1y_0^2y_2$ is also zero due to the corank $2$ condition.   
From now on, the equation of $S$ can be written as
\[
\begin{split}
    x_1^2 [c_0(y_1,y_2) + y_0 q(y_1,y_2)+ \mu y_0^2y_2] + x_0 x_1 [c_1(y_1,y_2) + 
  y_0 q_1(y_1,y_2)] \\+ x_0^2 y_2 (q(y_1,y_2) + \nu y_0y_2)=0.
\end{split}
\]
Moreover, if the coefficient of $x_0x_1y_0y_1^2$ is nonzero, then the singularity of $p$ has weight no worse than $(\frac{1}{2},\frac{1}{4},\frac{3}{8})$, which is impossible.   

\end{proof}

\subsection{Geometry of unstable surface}

\begin{theorem} \label{thm:unstable}
A bidegree-$(2,3)$-hypersurface $S\subseteq \PP^1\times \PP^2$ is unstable if and only if one of the following conditions holds 
\begin{itemize}
       \item[\rm(U1)]  $S=\PP^2\cup S'$ for some  $S' \in |\cO_{\PP^1\times \PP^2}(1,3)|$, meeting along a cuspidal cubic curve. 
      
      \item[\rm (U2)]   $S$ is singular along a line in the fiber whose projection under $\pi'$ is a triple line $3L$.
      And $\bB(S)= 2\ell \cup B'$ with $L\cap B' $ is a quartic point $\{o\}$.

      \item[\rm (U3)]   $S$ contains at least a corank $3$ isolated singularity $p$ and the fiber over $\pi(p)$ is a triple line $3L$. 
      And the projective tangent cone at $\pi^{\prime}(p)$ of branching locus $\PTC_{\pi^{\prime}(p)}\bB(S)$ contains 
      $3L$.

      \item[\rm (U4)]    $S$ contains at least a corank $3$ isolated singularity $p$ and the fiber over $\pi(p)$ is the union of a double line $2 L_1$ and a line $L_2$ meeting at $\pi'(p)$.
     And $\mult_{\pi'(p)}(\pi'(F),L_1)\geq 2$ for any fiber $F$.
     The projective tangent cone at $\pi^{\prime}(p)$ of the branching locus $\PTC_{\pi^{\prime}(p)}\bB(S)$ is at least $3L_1\cup L'$.
          
      \item[\rm (U5)] 
      $S$ is singular alone a type horizontal line $C$ with $\pi'(C)$ is a point $o\in \PP^2$.
      $S$ has a fiber which is the union of three lines intersecting at $o$.
      The projective tangent cone at $o$ of branching locus $\PTC_o\bB(S)$ is at least a quadruple line.

      \item [\rm (U6)] 
     $S$ is singular along a type horizontal line $C$ with $\pi'(C)$ is a point $o\in \PP^2$.
     $S$ has a  fiber $F_0$ such that the projective tangent cone of $\pi'(F_0)$ at $o$ is a double line $\PTC_o(\pi'(F_0))=2\ell$.
     And $\mult_o (\pi'(F),L)\geq 3$ for any fiber $F$.
     The projective cone $ \PTC_o\bB(S)$ is at least $3L\cup L'$.
     

      \item[\rm (U7)] 
      $S$ is singular along a type vertical line whose projection under $\pi'$ is a line $L$ .
      For any $F$ with $L\nsubseteq \pi'(F)$,  $\{o\} = L\cap \pi'(F)$ is at least a triple point.
      The branching locus $\bB(S) = 2\ell \cup B'$ with $L\cap B'$ containing at least a triple point which is exactly $o$ when such $F$ exists.

\end{itemize}
\end{theorem}

\begin{proof}
As in the proof of Theorem \ref{thm:semistable}, we present explicit proof for three slightly complicated cases.
\begin{enumerate}

    \item[(1)] If $S$ is of type (U4), then the equation of $S$ is given by
\[
  x_1^2 [c_0(y_1,y_2) + y_0 q(y_1,y_2) + \mu y_0^2 y_2 ] + x_0x_1 [c_1(y_1,y_2) + \nu y_0 y_2^2]  +x_0^2 y_2^2 l(y_1, y_2)=0.
\]
One observe that $S$ contains a corank $3$ isolated singularity $p\colon(1,0,1,0,0)$, 
and the fiber over $\pi(p)$ is the union of a double line $2 L_1\colon\{y_2^2=0\}$ and a line $L_2 \colon \{ l(y_1,y_2) =0\} $ meeting at $\pi'(p)$. 
One can assume that the coefficient of 
$x_1^2y_0^2y_2$ is nonzero, and the coefficients of $x_0^2y_2^3$, $x_0^2y_1y_2^2$ are not simultaneously zero.
Otherwise, the equation will degenerate to type (U5) and (U1). 
Then $\mult_{(1,0,0)} (\pi'(F),y_2)\geq 2$ for any fiber $F$ and $$\PTC_{\pi'(p)}\bB(S) = \nu^2 y_2^4-\mu y_2^3 l(y_1,y_2)= 3L_1 \cup L'.$$

Conversely, suppose that $S$ contains at least a corank $3$ isolated singularity $p$  and the fiber containing $p$ is a union of a double line $2 L_1$ and a line $L_2$ meeting at $p$. 
One can assume that $p= (1,0,1,0,0)$ and $L_1$ is given by $\{y_2=0\}$. 
It's easy to show that there are no $x_i x_j y_0^2 y_1$ terms in the equation of $S$ since 
 $\mult_{(1,0,0)}(\pi'(F),y_2)\geq 2$ for any fiber $F$. 
Then the equation of $S$ can be written as
\[
    x_1^2 [c_0(y_1,y_2) + y_0 q(y_1,y_2) + \mu y_0^2 y_2] + x_0x_1 [c_1(y_1,y_2) +y_0 q_1(y_1,y_2)] +x_0^2 y_2^2 l_0(y_1,y_2)=0.
\]   
The $\bB(S)$ is given by     
$$(c_1(y_1,y_2)+y_0 q_1(y_1,y_2))^2 - y_2^2 l_0(y_1,y_2)(c_0(y_1,y_2)+y_0 q_0(y_1,y_2) + \mu y_0^2 y_2)=0.$$   
Thus one can deduce that the $q_1= \nu y_2^2$ since $\PTC_{(1,0,0)}\bB(S)= 3L_1 \cup L'$.

    \item[(2)]  If $S$ is of type (U5), then the equation of $S$ is given by 
\[
   x_1^2 [c_0(y_1,y_2) + y_0 q_1(y_1,y_2)]
    + x_0 x_1 [c_1(y_1,y_2)+ \mu y_0 y_2^2]
    +x_0^2 c_2(y_1,y_2)=0.
\]
   We see at once that $S$ is singular along a horizontal line 
   $C \colon \PP^1\times (1,0,0)$ with $\pi'(C)$ is a point $o \colon (1,0,0) \in \PP^2$. 
   And the fiber over $(1,0) \in \PP^1$ is given by $\{c_2(y_1,y_2)=0\}$, which is the union of three lines intersecting at $o$. 
   In addition, $\PTC_o\bB(S)$ is a quadruple line $\{y_2^4=0\}$ if $\mu$ is nonzero. 
   Otherwise, the degree of $\PTC_o\bB(S)$ is at least $5$.

 Conversely, up to a coordinate change, one can assume that $S$ is singular along a horizontal line $C\colon\PP^1\times (1,0,0)$ with $\pi'(C)$ is a point $o \colon (1,0,0) \in \PP^2$. 
 In addition, the fiber over $(1,0) \in \PP^1$ is the union of three lines intersecting at $o$.  
 So the equation of $S$ can be written as
\[
   x_1^2 [c_0(y_1,y_2) + y_0 q_1(y_1,y_2)]
    + x_0 x_1 [c_1(y_1,y_2)+ y_0 q_2(y_1,y_2)]
    +x_0^2 c_2(y_1,y_2)=0.
\]   
 If $q_2 = 0$, then obviously $S$ is of type (U5). If $q_2 \neq 0$, then $\PTC_o\bB(S)$ is given by $q_2(y_1,y_2)^2$. Thus one can assume $q_2(y_1,y_2)= \mu y_2^2$ up to a coordinate change of $y_1$ and $y_2$ since $\PTC_o\bB(S)$ is a quadruple line.
    
  
    
    \item[(3)]  If $S$ is of type (U7), then the equation of $S$ is given by 
\[
\begin{split}
     x_1^2[c(y_1,y_2)+ y_0y_2 l(y_1,y_2) + \mu y_0^2y_2] + x_0x_1 (y_2^2l_0(y_0,y_1,y_2) + \nu y_1^2y_2 ) \\
    + x_0^2 y_2^2l_1(y_0,y_1,y_2)=0
\end{split}
\]
with $\bB(S)$  given by  
 \[
  y_2^2(y_2l_0(y_0,y_1,y_2)+\nu y_1^2)^2-y_2^2l_1(y_0,y_1,y_2)(c(y_1,y_2)+ y_0y_2 l(y_1,y_2) + \mu y_0^2y_2)=0.
\]
 One observes that $S$ is singular along a vertical line 
 $ \{x_1= y_2=0\}$ with $L\colon \{y_2=0\}$ the projection to $\PP^2$. 
 Also one can notice that $\bB(S) = 2L \cup B'$.
 If the coefficient of $x_1^2y_1^3$ is zero, then the projection  of every fiber under $\pi'$  contains $L$. 
 One can see that $L\cap B'$ is a quadruple point $(1,0,0)$.
 If the coefficient of $x_1^2y_1^3$ is nonzero, it is clear that $L\cap \pi'(F)$ is a triple point $o \colon (1,0,0)$ for any fiber $F$. 
 Finally, $L \cap B'$ is given by $\{\nu^2 y_1^4-y_1^3l_1(y_0,y_1)=y_2=0\}$, which contains a triple point $(1,0,0)$.

 Conversely, up to a coordinate change one may assume that $S$ is singular along a horizontal line $ \{x_1= y_2=0\}$ with $L\colon \{y_2=0\}$ the projection to $\PP^2$. 
 One can then write the equation of $S$ as
\[
\begin{split}
 x_1^2[c(y_0,y_1)+ y_2 q_0(y_0,y_1)+ y_2^2 l_0(y_0,y_1) + y_2^3 ] + x_0x_1 (y_2^2 l_1(y_0,y_1,y_2) + y_2 q(y_0,y_1) ) \\
    + x_0^2 y_2^2 l_2(y_0,y_1,y_2)=0.   
\end{split}
\]
 If the projection of every fiber  under $\pi'$ contains $L$,
 then there are no $x_1^2 c(y_0,y_1)$ terms.
 One can compute $\bB(S) = 2L \cup B'$ with $ B'\cap L = \{q(y_0,y_1)^2=0\}\subset \PP^1 $.
 By our assumption, we see $q(y_0,y_1) = l(y_0,y_1)^2$ and we may set $l = y_1$ by making coordinate change. This completes the proof.

 From here we assume that there is a fiber $F$ whose projection under $\pi'$ does not contain $L$. 
 One can assume that $F$ is over $(0,1)$. 
 Then $L\cap \pi'(F)$ is given by $\{c(y_0,y_1)=0\}$. 
 By our assumption, we see $c(y_0,y_1)=y_1^3$ and the triple point $o$ is $(1,0,0)$ after applying coordinate change.
Still $\bB(S) = 2L \cup B'$ and $L \cap B'$ is given by $\{q(y_0,y_1)^2- y_1^3 l_2(y_0,y_1,0)=0\}$.
It follows that $q=y_1^2$ since $L \cap B'$ contains the triple point $o$. 
\end{enumerate}

\end{proof}

\section{Minimal orbits}\label{section minimal orbit}

In this section, we study the strictly semistable locus of the GIT compactification .
According to Subsection \ref{numerical}, it suffices to discuss the points of type $(\alpha)-(\delta)$. 
The results are obtained by using the following criterion of Luna.
\begin{lemma}[Luna’s criterion \cite{Luna1975}]
 Let $G$ be a reductive group acting on an affine
variety $X$. If $H$ is a reductive subgroup of $G$ and $x \in X$ is stabilized by $H$, then the orbit $G \cdot x$ is closed if and only if $C_G(H) \cdot x$ is closed. \end{lemma}

As an immediate application, one can calculate the dimensions of the boundary strata (see subsection \ref{proof of thm:main}).

\begin{lemma}
 A generic cubic fourfold of type $(\alpha) - (\delta)$ gives a closed orbit. Therefore, each of the types $(\alpha) - (\delta)$ gives an irreducible boundary component for the GIT compactification.  The dimensions of these boundary strata are $4$ for $(\alpha)$, $3$ for $(\eta)$, $2$ for $(\gamma)$, and $1$ for $(\beta)$ and $(\delta)$.
\end{lemma}

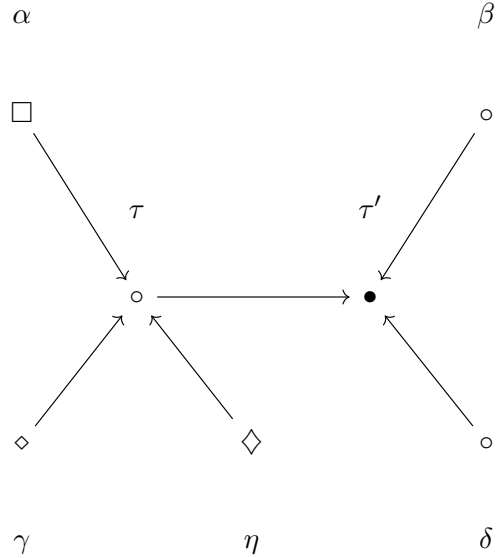
\begin{figure}
    \centering
\begin{tikzcd}
\alpha               &                  &                    &               & \beta             \\
\square \arrow[rdd]  &                  &                    &               & \circ \arrow[ldd] \\
                     & \tau             &                    & \tau^{\prime} &                   \\
                     & \circ \arrow[rr] &                    & \bullet       &                   \\
                     &                  &                    &               &                   \\
\diamond \arrow[ruu]  &                  &\diamondsuit \arrow[luu]&               & \circ \arrow[luu] \\
\gamma               &                  &\eta                &               & \delta             
\end{tikzcd}
    \caption{Incidence of the boundary components}
    \label{fig:my_label}
\end{figure}

Before starting the detailed analysis of the strata $(\alpha) - (\delta)$, we 
 describe the relations of their common degeneration. One observes that the types $(\alpha)$, $(\gamma)$ and $(\eta)$ have a common specialization  $$(\tau) \colon a_1x_1^2 y_0^2y_1+ a_2x_0x_1y_0y_1y_2 + a_3x_0^2y_1y_2^2 =0.$$
The stabilizer of $(\tau)$ contains a 1-PS $(4, -4, 3, 2, -5)$. Its equation further degenerates (for $a_1a_3 = 0$) to  $$(\tau^{\prime}) \colon x_0 x_1 y_0 y_1 y_2=0. $$
In addition, $(\tau^{\prime})$ is also a specialization of
the cases $(\beta)$, $(\gamma)$ and $(\delta)$. The resulting incidence diagram is given in the figure \ref{fig:my_label}.

\begin{lemma}
  A generic hypersurface of type $(\tau) $ and $(\tau^{\prime})$ is semi-stable with closed orbit.
\end{lemma}
\begin{proof}
This follows from Luna’s criterion cited above. 
The stabilizer of $(\tau)$ and $(\tau^{\prime})$ both contain a 1-PS $H$ of distinct weights. 
Thus it suffices to check the semi-stability with respect to the standard maximal torus $T= C_G(H)$ in $G$. 
The proposition follows easily.  
For example, the fact that $(\tau)$ is semi-stable follows from the simple observation that $(-2a+2b+c) + (2a+c+2(-b-c) )=0$ implies  either $-2a+2b+c \geq 0$ or $2a+c+2(-b-c) \geq 0$, 
where $(a,-a,b,c,-(b+c))$ are the weights of a 1-PS of $T$.
\end{proof}

We now do the case-by-case analysis of the minimal orbits of type $(\alpha)-(\delta)$. The common feature of all these cases is that the analysis reduces to some lower dimensional GIT problem.
Generally speaking, as the dimension of stratum increases, the analysis gets more involved. 
We first consider several cases with higher dimension.

   \begin{proposition}\label{prop1} 
Let $S$ be a hypersurface of type $(\alpha)$. 
Then $S$ is singular along a vertical line $L$ in fiber $F$.
And $S$ has a corank $3$ singularity $p$ such that $p \notin F$ and $\pi^{\prime}(p) \notin \pi'(L)$.
Moreover, we have
\begin{enumerate}  
    \item[(1)] $S$ is unstable if one of the following conditions hold
    \begin{enumerate}
        \item  $S=\PP^2\cup S'$ for some  $S' \in |\cO_{\PP^1\times \PP^2}(1,3)|$, where $\PP^2 = \pi^{-1}(\pi(F))$. And $\PP^2\cap S' = l_1^2 \pi'(L)$, where $l_1$ is a line in $\PP^2$.  
       
        \item
        $S$ has an another corank $3$ singularity $o$ such that $o \in L$.
        And $S$ is singular along a horizontal line containing $o$. In addition, the branch locus $\bB(S)$ contains a triple line.
        \item The fiber over $\pi(p)$ is a triple line $3L_1$ and the branch locus $\bB(S)$ also contains $3L_1$.
    \end{enumerate}
    \item[(2)] The orbit of $S$ is not closed if one of the following conditions hold
    \begin{enumerate}
        \item $S$ is singular along a vertical line in the fiber over $\pi(p)$. It degenerates to type $(\tau)$.
        \item $S=(\PP^1\times \PP^1) \cup S'$ for some  $S' \in |\cO_{\PP^1\times \PP^2}(2,2)|$. It degenerates to type $(\tau)$.
        \item $S=\PP^2\cup S'$ for some  $S' \in |\cO_{\PP^1\times \PP^2}(1,3)|$. It degenerates to type $(\tau^{\prime})$.
    \end{enumerate}
\end{enumerate}
Otherwise, $S$ is semistable with closed orbit. 
\end{proposition}
 
\begin{proof}
By inspecting the equation of $(\alpha)$, it is easy to see that $S$ is singular along a vertical line 
$\{x_1=y_2=0\}$, whose projection under $\pi'$ is given by
$L\colon \{y_2=0\}$. And  $S$ has a corank $3$ singularity at $(0,1,0,0,1)$. 
The converse is trivial.

Since the stabilizer of type $(\alpha)$ contains a 1-PS:
 $$H_1=\{ \diag(t^{3},t^{-3},t^2,t^2,t^{-4}) \;|\;t\in \CC^{\ast}   \}$$
and the center $C_G(H_1)\cong \CC^{\ast} \times \GL_2(\CC)$,
 we can reduce our problem to a simpler GIT problem $V^{H_1}\q C_G(H_1)$ by Luna's criterion where $V^{H_1} = \langle x_0^u x_1^{2-u} y_0^v y_1^{3-u-v} y_2^u \rangle$. 
 Then following the Hilbert-Mumford criterion, by diagonalizing 1-PS $\lambda$ into the form
 $$\lambda(t)=\diag (t^a,t^{-a},t^{b},t^{c},t^{-b-c})  \ \ (b\geq c)$$
 we can deduce that $S$ is unstable if and only if its defining equation is one of the following forms up to a coordinate change of $y_0$ and $y_1$
 \begin{itemize}
     \item $x_1^2 c(y_0,y_1)+ \mu x_0x_1y_1^2y_2 = 0$.
      \item $x_1^2y_1^2 l(y_0,y_1) + \mu x_0x_1y_1^2y_2 + \nu x_0^2y_1y_2^2 = 0$.
       \item $\mu x_1^2y_1^3 + \nu x_0x_1y_1^2y_2 + x_0^2y_2^2 l(y_0,y_1)=0$.
 \end{itemize}
And the orbit of $S$ is not closed if its equation has one of the following forms
\begin{itemize}
    \item $x_1 (x_1 c(y_0,y_1)+x_0q(y_0,y_1) y_2)=0$. 
    \item $x_0( x_1 q(y_0,y_1) y_2+ x_0 l(y_0,y_1) y_2^2)=0$.
    \item $y_1(x_1^2 q_1(y_0,y_1)+x_0x_1 l_1(y_0,y_1) y_2
    + \mu x_0^2 y_2^2)=0$.
    \item $x_1^2 y_1^2 l_0(y_0,y_1)+ x_0x_1 y_1 y_2 l_1(y_0,y_1) + x_0^2 l_2(y_0,y_1) y_2^2=0$.

\end{itemize}
Then the first two cases degenerate to type $(\tau')$, the last two cases degenerates to type $(\tau)$.

The proof of the geometric description is similar to Theorem \ref{thm:semistable} and Theorem \ref{thm:unstable}. To give an example, we show the third unstable case.
If the equation of $S$ is given by $ \mu x_1^2y_1^3 + \nu x_0x_1y_1^2y_2 + x_0^2y_2^2 l(y_0,y_1)=0.$
Then it is clear that the fiber over $(0,1)$ is a triple line $\{y_1^3=0\}$ and the branch locus also contains $\{y_1^3=0\}$.

Conversely, if $S$ is of type $(\alpha)$ whose fiber over $\pi(p)$ is a triple line $3L_1$. Then one can assume that $L_1 \colon \{y_1 =0\}$ up to a coordinate change. Its equation has the form
$$x_1^2 y_1^3 +x_0x_1 q(y_0,y_1) y_2+ x_0^2 l(y_0,y_1) y_2^2=0$$
with $\bB(S)$  given by  
$ y_2^2(q(y_0,y_1)^2-y_1^3l(y_0,y_1))=0.$
Thus one can deduce that $q=y_1^2$ by our assumption.   In conclusion, the equation of $S$ can be written as 
$$x_1^2y_1^3 + x_0x_1y_1^2y_2 + x_0^2y_2^2 l(y_0,y_1)=0.$$
\end{proof}

\begin{proposition}
Let $S$ be a hypersurface of type $(\gamma)$, then $S=(\PP^1\times \PP^1) \cup S'$ for some  $S' \in |\cO_{\PP^1\times \PP^2}(2,2)|$ and $S$ is singular along a horizontal line $C$ such that the intersection of $\PP^1\times \PP^1$ and $C$ is empty. Moreover, we have

\begin{enumerate}
    \item $S$ is unstable if the intersection $S' \cap (\PP^1 \times \PP^1) $  is of type $\cO_{\PP^1\times \PP^2}(1,0)\cup \cO_{\PP^1\times \PP^2}(1,2)$ or $\cO_{\PP^1\times \PP^2}(0,1)\cup \cO_{\PP^1\times \PP^2}(2,1)$, where the two components meet at a double point.
  
    \item The orbit of $S$ is not closed if the intersection $S' \cap (\PP^1 \times \PP^1) $ is a singular curve. 
    It degenerates to type $(\tau)$.
  
\end{enumerate}  

Otherwise, $S$ is semistable with closed orbit.  
  \end{proposition}
\begin{proof}
  The stabilizer of type $(\gamma)$ contains a 1-PS:
 $$H_2=\{ \diag(1,1,t^1,t^1,t^{-2}) \;|\;t\in \CC^{\ast}   \}.$$
The center $C_G(H_2)\cong \SL_2(\CC) \times \SL_2(\CC) \times \CC^{\ast}$.
 By Luna’s criterion, we can reduce our problem to a simpler GIT analysis $V^{H_2}\q C_G(H_2)$ where $V^{H_2} = \langle x_0^u x_1^{2-u} y_0^v y_1^{2-v} y_2 \rangle$. 
 Any 1-PS $\lambda$ can be diagonalized in the form
 $$\lambda(t)=\diag (t^a,t^{-a},t^{b},t^{c},t^{-b-c}).$$
 where $a\geq 0$ and $b\geq c$. Then following the Hilbert-Mumford criterion, we deduce that $S$ is unstable if and only if its defining equation is one of the following forms after we make a coordinate change.
 \begin{itemize}

     \item $y_1 y_2 (x_1^2 l(y_0,y_1)+ \mu x_0x_1y_1 + \nu  x_0^2 y_1) = 0$.
      \item $x_1y_2 (x_1 q(y_0,y_1)+ \mu x_0 y_1^2)= 0$.
 \end{itemize}
And the orbit of $S$ is not closed if its equation has the form
\begin{itemize}
    \item $y_2 (x_1^2 q(y_0,y_1) +x_0x_1 y_1 l(y_0,y_1) +\mu  x_0^2 y_1^2)=0$.
    \item $y_1 y_2 (x_1^2 l_0(y_0,y_1)+ x_0x_1l_1(y_0,y_1) + x_0^2 l_2(y_0,y_1)) = 0$.
    \item $x_1y_2 (x_1 q_1(y_0,y_1)+x_0 q_0(y_0,y_1))= 0$.
\end{itemize}
Then the first case degenerates to type $(\tau)$, and the other two cases are unstable. 

The remaining part follows by an analysis of the geometric description of the possible degenerations similar to before.
\end{proof}  
  
The remaining cases are quite similar, we omit the details.

\begin{proposition}
Let $S$ be a hypersurface of type $(\eta)$. Then $S$ has two isolated $\widetilde{E}_7$-type singularities $p$ and $q$ such that $\pi(p) \neq \pi(q)$ and $\pi^{\prime}(p) \neq \pi^{\prime}(q)$. The fibers over $\pi(p)$ and $\pi(q)$ contain a line $L_1$ and $L_2$ respectively with $L_1 \neq L_2$. And $\mult_{\pi'(p)} (\pi'(F),L_1)\geq 2$,  $\mult_{\pi'(q)} (\pi'(F),L_2)\geq 2$ for any fiber $F$. Moreover, we have   
   \begin{enumerate}
    \item[(1)] $S$ is unstable if one of the following conditions holds
    \begin{enumerate}
        \item $S=\PP^2\cup S'$ for some  $S' \in |\cO_{\PP^1\times \PP^2}(1,3)|$, meeting along a triple line.
        \item The fibers over $p$ and $q$ contain double line $2L_1$ and $2L_2$ respectively, with $\bB(S) = 3L_1 \cup 3L_2$. 
        \item The fibers over $p$ and $q$ are $l^2 L_1$ and $l^2L_2$ respectively, where $l$ is a line. And the branch locus $\bB(S)$ contains $4l$.
       
    \end{enumerate}
  
   \item[(2)] The orbit of $S$ is not closed if the fibers over $p$ or $q$ contains a double line $2l$ and the branch locus also contains the double line $2l$. It degenerates to type $(\tau)$.
  
\end{enumerate}  
 Otherwise, it is semistable with closed orbit. 
   
  \end{proposition}

\begin{proposition}
Let $S$ be a hypersurface of type $(\beta)$, then $S$ has two isolated $\widetilde{E}_8$-type singularities $p$ and $q$ such that $\pi(p) \neq \pi(q)$ and $\pi^{\prime}(p) \neq \pi^{\prime}(q)$.
The fibers over $p$ and $q$ are $3 L_1$ and $3 L_2$ respectively with $L_1 \neq L_2$.
Moreover, we have 
\begin{enumerate}
    \item $S$ is unstable if the branch locus contains a triple line.
  
    \item The orbit of $S$ is not closed if the branch locus contains a double line. It degenerates to type $(\tau^{\prime})$.
\end{enumerate}  
Otherwise, $S$ is semistable with closed orbit. 
\end{proposition}

\begin{proposition}\label{prop5}
 Let $S$ be a hypersurface of type $(\delta)$, then $S= \PP^2 \cup \PP^2 \cup (\PP^1 \times C)$ for some cubic curves $C$ in $\PP^2$. Moreover, we have 
\begin{enumerate}
    \item $S$ is unstable if  if $C$ has cusps or triple points, or is reducible with two components tangent at a point.

    \item The orbit of $S$ is not closed if $C$ is the union of a conic and a transversal line.
\end{enumerate} 
   Otherwise, it's semistable with closed orbit.

\end{proposition}

  \begin{proof}
 One can deduce this from the well-known GIT analysis of plane cubic curves, for example, see \cite[Section 1.11]{Mumford1977}.
  
\end{proof}
  
\section{Boundary of \texorpdfstring{$\overline{\cM}_2^{\rm GIT}$}{M2GIT}}\label{section boundary}

\subsection{Stable surface with isolated singularity}
The following is a direct consequence of Theorem \ref{thm:semistable}. 
\begin{proposition}\label{isolated}
 Let  $S\in |\cO_{\PP^1\times \PP^2}(2,3)|$ with only isolated singularity. 
 Then $S$ with non-\rm ADE singularities is stable if and only if one of the following situations holds
  \begin{enumerate}
     \item [$(\zeta)$] $S$ has a $\widetilde{E}_7$-type singular point $p$ such that there exists a fiber $F$ with $\pi'(p) \notin \pi'(F)$. 
     \item [$(\xi)$] $S$ has a  $\widetilde{E}_8$-type singular point $p$ such that the fiber containing $p$ is not a triple line. 
 \end{enumerate}
\end{proposition} 

\begin{proof}
Suppose $S$ has only at worst isolated singularities of type $(\zeta)$ or $(\eta)$. By Theorem \ref{thm:semistable}, we immediately know that $S$ is stable. 
Conversely, suppose $S$ is stable and it has a non-ADE corank $2$ isolated singularity at $p\colon(1,0,1,0,0)$,
then $p$ is an isolated $\widetilde{E}_7$ or $\widetilde{E}_8$ singularity. 
Up to a change of coordinates preserving $p$, the local defining equation near $p$ has three possibilities :
\begin{itemize}
    \item $\mu x_1^2 + 3$ jets $+ 4$ jets $+ 5$ jets,
    \item $\mu y_1^2 + 3$ jets $+ 4$ jets $+ 5$ jets,
    \item $\mu(x_1+y_1)^2 + 3$ jets $+ 4$ jets $+ 5$ jets.
\end{itemize}
Direct computation shows the only possible case is the third one.

   If $p$ is $\widetilde{E}_7$ type, take the change of coordinate: $x_1+y_1 \mapsto z$. Then there are no terms $y_1^k y_2^{3-k}$ in the third jets. The equation of $S$ is of the form 
 \begin{equation}\label{eq E7 tilde }
   \begin{aligned}
       f= x_1^2(c_0(y_1,y_2) + y_0 q_0(y_1,y_2)+y_0^2 l_0(y_1,y_2)) + x_0x_1(c_1(y_1,y_2) + y_0 q_1(y_1,y_2))\\ + x_0^2(y_1^2 l_1(y_1,y_2) + y_1y_2^2) +\mu(x_1^2 y_0^3 + 2 x_0 x_1 y_0^2 y_1 + x_0^2 y_0y_1^2), 
   \end{aligned}  
 \end{equation}     
     satisfying  $a_{210} - a_{120} + a_{030} = 0$, $a_{201} - a_{111} + a_{021} = 0 $ and $a_{102} = a_{012}$, where $a_{ijk}$ represents the coefficient of $x_0^{2-i}x_1^i y_0^{3-j-k} y_1^j y_2^k$.
    Conversely, for the general equations of the form (\ref{eq E7 tilde }), take the change of coordinate: $x_1+y_1 \mapsto z$. Then the weights on variables $z$, $y_1$ and $y_2$ are $(\frac{1}{2},\frac{1}{4},\frac{1}{4})$.
Thus $p$ is $\widetilde{E}_7$ type. 
And one can see that $\pi'(F)$ does not contain $\pi'(p)=(1,0,0)$, where $F$ is the fiber over $(0,1)$.

If $p$ is $\widetilde{E}_8$ type, take the change of coordinate: $x_1+y_1 \mapsto z$. Suppose that the weight of $y_1$ is less than $\frac{1}{6}$, then there are no terms $y_1^3$, $y_2^2 y_1$, $y_2 y_1^2$, $z y_1^2$ in the third jets.
No term $y_1^4$, $y_1^3 y_2$, and $y_1^5$ in the fourth and fifth jets. Thus the equation of $S$ is of the form 
\begin{equation}\label{eq E8 tilde}
 \begin{aligned}
   f= x_1^2(y_2^3 + y_2^2y_1 +y_2y_1^2 + y_0 q_0(y_1,y_2)+y_0^2 l_0(y_1,y_2) )\\ + x_0x_1(c_0(y_1,y_2) + y_0 q_1(y_1,y_2) ) + x_0^2c_1(y_1,y_2) \\ 
     + \mu(x_1^2 y_0^3 + 2 x_0 x_1 y_0^2 y_1 + x_0^2 y_0y_1^2),
 \end{aligned}   
\end{equation}
     satisfying $a_{201} - a_{111} + a_{021} = 0 $, $a_{220} = a_{130}$, $a_{211} = a_{121}$, 
     $a_{120} = 2 a_{210} = 2a_{030}$ and $a_{102} = a_{012}$.   
Conversely, for the general equations of the form (\ref{eq E8 tilde}), take the change of coordinate: $x_1+y_1 \mapsto z$. Then the weights on variables $z$, $y_1$ and $y_2$ are $(\frac{1}{2},\frac{1}{6},\frac{1}{3})$. And one can see that the fiber containing $p$ is not a triple line.

Similarly, if the weight of $y_2$ is less than $\frac{1}{6}$, then equation of $S$ is of the form 
\[\begin{split}
   f= x_1^2(c_0(y_1,y_2) + y_0 q_0(y_1,y_2)+y_0^2 l_0(y_1,y_2) ) +  x_0^2 y_1^2 l_1(y_1,y_2)\\ +x_0x_1(y_1^3 + y_1^2y_2 +y_1y_2^2 + y_0 y_1 l_2 (y_1,y_2) ) 
     \\ + \mu(x_1^2 y_0^3 + 2 x_0 x_1 y_0^2 y_1 + x_0^2 y_0y_1^2).
     \end{split}
     \]
     satisfying $a_{201} - a_{111} + a_{021} = 0 $. One can easily see that 
     $S$ is singular along $\{x_1=y_1=0\}$. Thus $p$ is not an isolated singularity.
\end{proof}

\subsection{Stable surface with non-isolated singularity}

\begin{proposition}\label{non-isolated}
 Let $S$ be a bidegree $(2,3)$-hypersurface of $\PP^1\times \PP^2$ with non-isolated singularities. Then one of the following situations holds: 
 \begin{enumerate}
    \item $\mathrm{Sing}(S)$ contains a vertical line and $S$ is not properly stable. 
    \item $\mathrm{Sing}(S)$ contains a horizontal line and $S$ is not properly stable. 
    \item $\mathrm{Sing}(S)$ contains an elliptic curve and $S$ is not properly stable unless $S$ is of type $(1,1)+(1,2)$ or $(0,2)+(2,1)$.
    \item $\mathrm{Sing}(S)$ contains a section of degree $1$ and $S$ is stable.  
    \item $\mathrm{Sing}(S)$ contains a section of degree $2$ and $S$ is stable. 
     \end{enumerate}
\end{proposition}

\begin{proof}

Let $C\subseteq \mathrm{Sing}(S)$ be an irreducible curve.  If $S=S_1\cup S_2$ is reducible, one can assume the bidegree of $S_1$ is $(1,0)$, $(0,1)$, $(0,2)$ or $(1,1)$. This gives the elliptic curve cases. For the first two cases, $S$ is not properly stable by Theorem \ref{thm:semistable}.

If $S$ is irreducible, then $C$ is either a section or contained in a fiber as the intersection of  $C$ with any fiber is at most $1$ by genus calculation. Suppose that $C$ is contained in a fiber. Then the only possible case is that $S$ is singular along a vertical line by some simple computation. Next, we assume that $C$ is a section. Let $d$ be the degree of $C$ as a section. Take a general hypersurface $H$ of type $(0,1)$, then $H\cap S$ is an irreducible curve and it is singular along $H\cap C$. Since the arithmetic genus of $H\cap S$ is at most $2$, then $H\cap C$ has at most $2$ points. It follows that $d$ is at most $2$.

For $d=0$, the surface is just $\rm (N4)$ type, which is not properly stable.

For $d=1$, we claim that up to change of coordinate, the equation of $S$ is of the form
      \[
      \begin{split}
          (\theta)\colon  f=  y_2^2 c_{2,1}(x_0,x_1 | y_0,y_1,y_2) 
     +y_2 (x_0y_1 - x_1y_0) q_{1,1}(x_0,x_1 | y_0,y_1) \\
     +(x_0y_1 - x_1y_0)^2 l(y_0,y_1) ,
      \end{split}
       \]  
where $c_{2,1}(x_0,x_1 | y_0,y_1,y_2)$ represents the bidegree $(2,1)$ homogeneous polynomial of $x_0$, $x_1$ and $y_0$, $y_1$, $y_2$.  
And $q_{1,1}(x_0,x_1\, |\, y_0,y_1)$ is similar. 
After changing coordinates, we may assume that $S$ is singular along a curve $C$ of the form $\{ (x_0,x_1,x_0,x_1,0) | (x_0,x_1)\in \PP^1\}$. 
Its equation is given by  
$C \colon\{ y_2=x_0y_1-x_1y_0=0\}.$  Since $S$ contains $C$, the equation of $S$
is of the form $f = y_2 g + (x_0y_1 - x_1y_0) h $ where $g \in  |\cO_{\PP^1\times \PP^2}(2,2)|$ and $h \in |\cO_{\PP^1\times \PP^1}(1,2)|$. Via computing $\frac{\partial f}{\partial y_2}$ and $\frac{\partial f}{\partial y_1}$ on $C$, one can deduce that both $\{ g=0 \} $ and $\{ x_0 h = 0 \}$ contain $C$. This gives the equation of type $(\theta)$.
Conversely, if the equation of $S$ is of type $(\theta)$, then the Jacobi of $S$ is zero along $C$.

For $d=2$, we claim that the general equation is of the form
      \[\begin{split}
         (\phi)\colon f= (x_0y_1 - x_1y_0)^2 l_0(y_0,y_1,y_2) + (x_0y_1 - x_1y_0)(x_0y_2-x_1y_1) l_1(y_0,y_1,y_2)  \\+ (x_0y_2 - x_1y_1)^2 l_2(y_0,y_1,y_2). 
      \end{split}
     \]   
As above, one can assume that $S$ is singular along the curve $C \colon \left\{ x_0y_2 - x_1y_1= x_0y_1 - x_1y_0=0\right\}$.  
Since $C$ is contained in $S$, the equation of $S$ is of the form 
$f = (x_0y_1 - x_1y_0) g + (x_0y_2 - x_1y_1) h$, 
where $g$, $h \in  |\cO_{\PP^1\times \PP^2}(1,2)|$.
One can deduce that both $ \{x_1 g=0 \}$ and $\{x_0 h=0 \} $ contain $C$ because $\frac{\partial f}{\partial y_2}$ and $\frac{\partial f}{\partial y_0}$ are zero on $C$. 
Furthermore, using the equation of $C$, it turns out the only possibility is that both $\{ g=0\} $ and $\{h=0 \}$ contain $C$. 
Thus the equation of $S$ is of type $(\phi)$. 
Conversely, if the equation of $S$ is of type $(\phi)$, the Jacobi of $S$ is zero along $C$.
Finally, the assertion of stability follows directly from Theorem \ref{thm:semistable}.
\end{proof}

\subsection{Proof of Theorem \ref{thm:main}}\label{proof of thm:main}
Note that the normal surface singularity admits  crepant resolution if and only if it is of ADE type.
So $\cM_2$ consists of $(2,3)$-hypersurfaces with only simple singularities and by Theorem \ref{thm:semistable}, we know that $\cM_2$ is lying in the stable locus of $\overline{\cM}_2^{\rm GIT}$.
The description of the boundary components simply follows from the combination of Proposition \ref{prop1}-\ref{prop5}, \ref{isolated} and \ref{non-isolated}.
The dimension of the strictly semistable boundary $(\alpha)-(\delta)$ can be computed via Luna’s criterion. For example, $$\dim  (\alpha)= \dim \mathbb{P}V^{H_1} \q  C_G(H_1) =4, $$
where $V^{H_1} = \langle x_0^u x_1^{2-u} y_0^v y_1^{3-u-v} y_2^u \rangle$ and $C_G(H_1)\cong \CC^{\ast}  \times \GL_2(\CC)$. For stable components, we can also compute the dimension as follows:
\begin{itemize}
    \item For $(\zeta)$, it can be viewed as the quotient space $\PP(V_1)/G_1$, where $V_1$ is the vector space
spanned by monomials in the equation $(\zeta)$ and $G_1$ is the group fixing the singular point $p\colon (1,0,1,0,0)$ and the $2$-jets $(x_1+y_1)^2$. 
As $\dim \PP(V_1)=16 $ and $\dim G = 6$, we get $\dim (\zeta) = 10$. Similarly, $(\xi)$ is the quotient space $\PP(V_2)/G_2$ with $\dim \PP(V_2)=13 $ and $\dim G_2 = 6$. It follows that $\dim (\xi) = 7$.
    \item Similar as above, $(\theta)$ is the quotient space $\PP(V_3)/G_3$, where $V_3$ is the vector space
spanned by monomials in the equation $(\theta)$ and $G_3$ is the group fixing the section of degree $1$ with the parametric form $\{ (x_0,x_1,x_0,x_1,0) | (x_0,x_1)\in \PP^1\}$. One can calculate that $\dim \PP(V_3)=14 $ and $\dim G_3 = 6$. 
It follows that $\dim (\theta)=8$. Similarly, one can see that 
$\dim (\phi)=5$.

    \item $(r_1)$ is a union of a $(1,1)$-hypersurface and a $(1,3)$-hypersurface. 
    The elements of $(r_1)$ are parametrized by the product of two projective spaces $\PP(V) \times \PP(V')$, where $V= H^0(\cO_{\PP^1 \times \PP^2}(1,1))$ and $V' = H^0( \cO_{\PP^1 \times \PP^2}(1,3))$. 
    Thus $$\dim (r_1)= \dim \PP(V)+\dim\PP(V')-\dim \PGL_2 \times \PGL_3 = 13.$$ With the same method, one can get  $\dim (r_2)=2$.
\end{itemize}
The last statement follows from Proposition \ref{prop:projective model}.
\qed

\section{Arithmetic compactifications} \label{section Arithmetic compactifications}
\subsection{Baily--Borel compactifications of $\cF_{T_n}$}

Recall that the type II and III boundary components of the Baily--Borel compactifications of $\cF_{T_n} $ correspond one-to-one to the classes of rank $2$ and $ 1 $  isotropic sublattices of $\Sigma_n$ modulo $ \Gamma_n $. 

Let $I_2(\Sigma_n) = \bigcup_{e} I_{2,e}(\Sigma_n)$ be the set-theoretic partition of rank $2$ isotropic sublattices of $\Sigma_n$ where 
\[
 I_{2,e}(\Sigma_n) = \fbrace{J\in I_2(\Sigma_n)}{|H_J| = e}{\Big}.  
\]
Here $H_J = {(J^{\perp})}^\perp_{\Sigma^*_n}$ with both the orthogonals taken inside $ \Sigma_n^*$, is a isotropic subgroup of $A_{\Sigma_n}$ with the induced finite quadratic form. 
Therefore $e^2 \,\big|\, (9 = \lvert A_{\Sigma_n}\rvert) $ and $e=1,3$.

\begin{lemma}\label{lem:good basis}
For $ n=1,2,3$ and $e=1,3$ fixed, all $J^\perp/J$ have the signature $ (0,16) $ and the same discriminant form $(A_{n,e},q_{n,e})$ for $ J\in I_{2,e}(\Sigma_n) $, i.e., are in the same genus which we denote by  $\cG(n,e) $.
Moreover, for $ J \in I_{2, e}(\Sigma_n)$ and  $e=1$ or $e=3$ , there exists $\left\{v_1, \ldots, v_{20}\right\}$ a basis of the lattice $\Sigma_n$ such that $J=\left\langle v_1, v_2\right\rangle, J^{\perp}=\left\langle v_1, \ldots, v_{18}\right\rangle$ and the quadratic form in this basis is of the form
$$
Q=\left(\begin{array}{ccc}
0 & 0 & A \\
0 & B & 0 \\
{ }^t A & 0 & D
\end{array}\right)
$$
where
$$
A=\left(\begin{array}{cc}
0 & 1 \\
e & 0
\end{array}\right), D=\left(\begin{array}{cc}
2t& 0 \\
0 & 0
\end{array}\right), 
$$
$0\leq t < e$ is uniquely determined by $n$  and $B$ is any matrix representing the quadratic form on $J^\perp /J$.

\end{lemma}

\begin{proof}
  It's well-known that we have the isomorphism of finite quadratic form \[
    A_{J^\perp/J} \cong H_J^\perp /H_J
     \]
(c.f. \cite[Proposition 6.5]{camere2018some} and \cite[Proof of Lemma 5.1.3]{Scattone1987}).
For $e=1$, we have $A_{J^\perp/J} \cong H_J^\perp /H_J\cong A_{\Sigma_n}$.
For $e=3$, one can find that $H_J^\perp  = H_J$ for $n=1,2,3$ and we have $A_{J^\perp/J} = \{0\}$.
So $\cG(n,e)$ is well-defined.

For the second statement, the proof of \cite[lemma 5.2.1]{Scattone1987} works and we sketch the proof here for the convenience of readers.
Since $J \subset J^{\perp}$ are primitive sublattices of $\Sigma_n$, 
  we can choose a basis of $\Sigma_n$ in which 
 $J=\left\langle v_1, v_2\right\rangle, J^{\perp}=\left\langle v_1, \ldots, v_{20}\right\rangle$  The matrix $Q(v_i,v_j) $ will then have the form 
 $$
Q=\begin{pmatrix}
0 & 0 & A_{0}\\
0 & B & C_{0}\\
A_{0}^{t} & C_{0}^{t} & D_{0}
\end{pmatrix}
$$ where $B$ represents the bilinear form of $J^\perp/J$.
Recalling that $ J \in I_{2,e}(\Sigma_n) $,  $H_J \cong Z/eZ $, a direct application of elementary divisors produces matrices $U$,$Z \in \GL_2(\mathbb Z)$ such that $U^t A_0 Z= \left( \begin{smallmatrix}
0 & 1\\
e & 0
\end{smallmatrix} \right)$.
Therefore, the change of 
basis described by the matrix $g=\operatorname{diag}(U,\mathrm{Id}_{16},Z) \in \GL_{20}(\mathbb Z) $ 
transform $A_0$ into $A$, $C_0$ into $C_1$ and preserves $B$.

Next Scattone showed that there exist integral matrices $V$ and $Y$ such that $BY+VA+C_1=0$ if $\gcd(e,\det B) =1$,  which is satisfied since $e^2\cdot \det B =9$ and we always have $e$ or $\det B$ equals $1$.
By choosing $V$ and $Y$ as above, and applying change of basis
$
 g=\begin{pmatrix}
I & V^{t} & \\
 & I & Y\\
 &  & I
\end{pmatrix}
$
we put Q into the form
$
\begin{pmatrix}
0 & 0 & A\\
0 & B & 0\\
A^{t} & 0 & D_{2}
\end{pmatrix}
$.
Finally, by applying 
$
 g=\begin{pmatrix}
I & 0 & W\\
 & I & 0\\
 &  & I
\end{pmatrix}
$
and choose a appropriate $W$, we can put $D_2$ into the required form.
For $e=1$, $t$ can only be $1$.
\end{proof}

Let $N_H(J)\defeq \im\big(\Stab_H(J) \xrightarrow{r} \GL(J) \cong \GL_2(\ZZ)\big)$ be image of the stabilizer of $J$ under the action of the group $H< O(\Sigma_n)$ in $\GL_2(\ZZ)$.
\begin{lemma}   \label{lem:stabilizer}
  For $J\in I_{2,e}(\Sigma_n) $, we have:
$$
N_{O(\Sigma_n)}( J) \cong \bigg\{\begin{pmatrix}
a & be\\
c & d
\end{pmatrix} 
\in \GL_{2}( Z) \;\bigg|\;  a^{2} \equiv 1 \ (\bmod \ e) \bigg\} \ 
$$
and
$$
N_{O^{+}(\Sigma_n)}( J) \cong \bigg\{
\begin{pmatrix}
a & be\\
c & d
\end{pmatrix} 
\in \SL_{2}( Z) \;\bigg|\;  a^{2} \equiv 1 \ (\bmod \ e) \bigg\} . 
$$
\end{lemma}
 
\begin{proof}
 For the case of polarized K3 surface, this is \cite[Lemma 5.6.3, 5.6.6]{Scattone1987}).
 The proof for our case is basically the same.
 Let $g$ be a general element of  $\Stab_{O(\Sigma_n)}(J)$ of the form
$$ 
  g=\begin{pmatrix}
        U & V & UW \\
          & X & Y  \\
          &   & Z
    \end{pmatrix}.
$$
Take $Q$ as in the previous lemma.
The condition $g^{t} Qg=Q$ gives  
 $$ A=U^{t} AZ , B=X^{t} BX,
 C= X^{t} BY + V^t A Z,
 D=Z^{t} DZ+Y^{t} BY+W^{t} U^{t} AZ+Z^{t} A^{t} UW
 $$
 which implies  $U\in N_{O(\Sigma_n)}(J)$ is of the form $ \begin{pmatrix}
a & be\\
c & d
\end{pmatrix} \in \GL_{2}(\mathbb{Z})$ and 
$
D \equiv \begin{pmatrix}
2a^{2} t & \ast \\
\ast  & \ast 
\end{pmatrix} \; ( \bmod \ 2e)
   $.
Note that in our case we have $ \gcd(t,e)=1 $, which implies 
$a^{2} \equiv 1 \ (\bmod \ e) $ is a necessary condition.
   
Conversely, if $U \in \GL_2(\mathbb Z) $ is  an arbitrary matrix 
   of the form 
$
    U =  \begin{pmatrix}
            a & be  \\
            c & d
        \end{pmatrix} 
    \in \GL_{2}(\mathbb{Z})  
$
   with  $a^{2} \equiv 1 \ (\bmod \ e) $, let $Z=A^{-1} \left( U^{-1}\right)^{t}  A $
 and then
 $$
   g=\begin{pmatrix}
 U & 0 & UW\\
 & I & 0\\
 &  & Z
    \end{pmatrix}
      $$
  gives a lift of $U$
   in $\Stab_{O( \Sigma_n)}(J)$ 
   with an appropriately chosen $W$.
   
   The second statement holds naturally by considering the sign of the determinant.

 \end{proof}

\begin{proposition}
 \label{prop:bij fiber}
  There is a bijection between
  \[
  I_{2,e}(\Sigma_n)\big\slash \Gamma_n \xlongleftrightarrow{1:1} \cG(n,e).
  \]
\end{proposition}

\begin{proof}
We divide the proof into three steps.

Firstly we show that 
 $ I_{2,e}(\Sigma_n)\big\slash O(\Sigma_n)\xlongleftrightarrow{1:1} \cG(n,e)$.  
  Let $L(e,t)$ be the lattice given by the Gram matrix 
  $\left(\begin{smallmatrix}
0 & e\\
e & 2t
\end{smallmatrix}\right)$. 
By Lemma \ref{lem:good basis} above, given $J\in I_{2,e}(\Sigma_n)$, one get 
 $ \Sigma_n \cong U   \oplus L(e,t) \oplus \left( J^{\perp }/J \right) $ where $J^\perp /J $ belongs to $\cG(n,e)$. 
 Conversely, for any $M \in  \cG(n,e)$
  since  $U  \oplus L(e,t) \oplus M $ and $\Sigma_n$ belong to the same
  genus, we have $ \Sigma_n \cong U   \oplus L(e,t) \oplus M$ so by the uniqueness of indefinite even lattice (see \cite[Corollary 1.13.3]{Nikulin1979}).
 With this isomorphism one can find $J\in I_{2,e}(\Sigma_n) $ such that $ J^{\perp }/J \cong M$.
 Therefore, 
  $$
  \cG(n,e) \rightarrow I_{2,e}(\Sigma_n)\big\slash O(\Sigma_n) 
  $$ is a bijection.

Next we show that \[
 I_{2,e}(\Sigma_n)\big\slash O^+(\Sigma_n) \rightarrow I_{2,e}(\Sigma_n)\big\slash O(\Sigma_n)
\] is a bijection.
Let $H_1 < H_2$ be a finite index subgroup acting on set $\mathscr{S}$.
Note that the fiber of $\mathscr{S}  /H_1 \rightarrow \mathscr{S}/H_2$ over $[a]/H_2$ has cardinality $\frac{[H_2:H_1]}{[\Stab_{H_2}(a):\Stab_{H_1}(a)]}$ for $a\in \mathscr{S}$.
In our case, $[O(\Sigma_n): O^+(\Sigma_n)] =2$.
By  lemma \ref{lem:stabilizer}
we have $N_{O^{+}(\Sigma_n)}( J) \neq  N_{O(\Sigma_n)}( J) $, hence $\Stab_{O^{+}(\Sigma_n)}( J) \neq  \Stab_{O(\Sigma_n)}( J)$.
Then one can conclude that
   $ I_2(\Sigma_n) / O(\Sigma_n) \cong I_2(\Sigma_n) / O^{+}(\Sigma_n) $.

Finally we consider the map  \[
I_{2,e}(\Sigma_n)\big\slash \Gamma_n \rightarrow I_{2,e}(\Sigma_n)\big\slash O^+(\Sigma_n).
\] 
For $n=1$ and $2$, one have $ O^{+}(\Sigma_n)/ \Gamma_n \cong O(A_{\Sigma_n}) = \{ \id,-\id \} $ .
Since $-\id$ doesn't change the $\Gamma_n$-class of $J$, the conclusion follows.

For $n=3$, we have $O^{+}(\Sigma_3)/ \Gamma_n  \cong O(\Sigma_3)/ \widetilde{O}(\Sigma_3) \cong O(A_{\Sigma_3})\cong \ZZ/2\ZZ \times \ZZ/2\ZZ $ by \cite[Theorem 1.14.2]{Nikulin1979}.
To determine $\Stab_{O^+}(J) /\Stab_{\Gamma_n}(J),$ we consider the exact sequence \[
1   \rightarrow \Fix_{O^+(\Sigma_n)}(J) /\Fix_{\Gamma_n}(J) \xrightarrow{i} \Stab_{O^+}(J) /\Stab_{\Gamma_n}(J) \xrightarrow{r} N_{O^+}(J) /N_{\Gamma_n} (J) \rightarrow 1 .
\]
Note that \[
\begin{tikzcd}
1 \arrow[r] & \ker \tilde{\beta} \arrow[r] \arrow[d, hook] & \Fix_{\Gamma_n}(J) \arrow[r, "\tilde{\beta}"] \arrow[d, hook] & \widetilde{O}^+(W_e) \arrow[d, hook] \arrow[r] & 1 \\
1 \arrow[r] & \ker \beta \arrow[r]                         & \Fix_{O^+(\Sigma_n)}(J) \arrow[r, "\beta"]                    & O^+(W_e) \arrow[r]                             & 1
\end{tikzcd}.
\]
We have \[
 1 \rightarrow  \ker \beta / \ker \tilde{\beta} \rightarrow \Fix_{O^+(\Sigma_n)}(J) /\Fix_{\Gamma_n}(J) \to O^+(W_e)/\widetilde{O}^+(W_e)
 \to 1
\]
for $W_e = J^\perp /J$ where $J\in I_{2,e}(\Sigma_n)$  and $\beta, \tilde{\beta}$ are the corresponding restriction maps.
  
For the case $e=1$, note that $\Sigma_3 \cong U^{ 2} \oplus W_1$, we have $\ker \beta  = \ker \tilde{\beta}$.
Hence $$\Fix_{O^+(\Sigma_n)}(J) /\Fix_{\Gamma_n}(J) \cong  O^+(W_e)/\widetilde{O}^+(W_e).
$$
Since there are only two types of even (negative) definite unimodular lattice of rank $16$, $D_{16} \subset D_{16}^{+}$ and $E_8 \oplus E_8$,
Since the only two isotropic subgroups of $A_{\Sigma_3}$ are differed by $\id \times (-\id)  \in O(A_{W_1}) $, one see that $\id \times (-\id)$ is in the image of $O(W_1) \rightarrow O(A_{W_1})$ by Lemma \ref{lem:overlattice} below, which generates $O(A_{W_1})$.
Hence $O(W_1) \rightarrow O(A_{W_1}) $ is surjective and $ O^+(W_1)/\widetilde{O}^+(W_1) \cong O(A_{\Sigma_3})$.
Form the exact sequence, the cardinality of the fiber over $[J]_{O^+(\Sigma_n)}$ is less or equal than 1 which gives the conclusion.

For $e=3$, the lattice $W_3$ is unimodular and we have $\Sigma_3 \cong U\oplus U(3) \oplus W_3$.
Using normalized basis, one can easily show that
$$
N_{\Gamma_3}( J) \cong \bigg\{
\begin{pmatrix}
a & 3b\\
c & d
\end{pmatrix} 
\in \SL_{2}( Z) \;\bigg|\;  a \equiv 1 \ (\bmod \ 3) \bigg\} . 
$$ and
$\big \lvert  N_{O^+}(J) /N_{\Gamma_3} (J) \big\rvert =2$ for $J\in I_{2,3}(\Sigma_3)$.
One the other hand we have \[
 \Fix_{O^+(\Sigma_n)}(J) /\Fix_{\Gamma_n}(J) \cong \big(\ker \beta / \ker \tilde \beta \big) \neq 1
\]
So $\big [ \Stab_{O^+}(J) : \Stab_{\Gamma_3}(J) \big ]  \geq 4$.
The conclusion follows from the cardinality formula of the fiber.

\end{proof}

 \begin{lemma}{\cite[Proposition 3.6]{Ebeling2013}}    \label{lem:overlattice}
    Let $\Lambda \subset \mathbb{R}^{n}$ be a  lattice. 
    There is a natural one-to-one correspondence between isomorphism classes of even overlattices $\Lambda \hookrightarrow \Gamma$  and orbits of isotropic subgroups $H \subset A_{\Lambda}$ under the image of the natural homomorphism $O(\Lambda) \rightarrow O\left(q_{\Lambda}\right)$. Unimodular lattices correspond to subgroups $H$ with $|H|^{2}=\left|A_{\Lambda}\right|$.
    
    In other words, there is a bijection between 
    isotropic subgroups of discriminant group $A_{\Lambda}$ and 
    overlattices of $\Lambda$.
  \end{lemma}

With Proposition \ref{prop:bij fiber} and the lemma, we have 
\begin{theorem}\label{thm:Baily--Borel}
Let $\cF_{T_n}^\ast$ be the Baily--Borel compactification of $\cF_{T_n}$. Then the boundary $\partial\cF_{T_n}:=\cF_{T_n}^\ast-\cF_{T_n}$ is given as follows:
\begin{enumerate}
    \item $n=1$, $\partial \cF_{T_1}$ consists of $14$ modular curves and $2$ points. 
    All curves meet at one point. There are $2$ curves meet at another point.
    
    \item $n=2$, $\partial \cF_{T_2}$ consists of  $11$ modular curves and $2$ points.
    All curves meet at one point. There are $2$ curves meet at another point.
    \item $n=3$, $\partial \cF_{T_3}$ consists of $10$  modular curves and $2$ points. 
    All curves meet at one point. There are $2$ curves meet at another point.
\end{enumerate}
\end{theorem}

\begin{proof}

By Proposition \ref{prop:bij fiber}, there is a one-to-one correspondence between isomorphic classes of lattices in  the genus $ \cG(n,e) $ and modular curves. 
To determine the classes $L$ in $\cG(n,e)$, by Nikulin's result \cite[Proposition 1.6.1]{Nikulin1979}, it's sufficient to classify all the the primitive embeddings of $F$ of signature $(0,8)$ with $A_L \cong A_F$ and $q_L \cong - q_F$, into some $H \in \uniL_{0,24}$ unimodular signature $(0,24)$ lattices up to $O(H)$-equivalence.

(i). For $n=1$, we have $(A_{\Sigma_1},q_{\Sigma_1}) \cong (A_{A_8}, -q_{A_8}) $.
 So we consider the all possible primitive embeddings of $A_8$ for $e=1$  and the primitive embedding of $E_8$ for $e=3$, into some $H\in \uniL_{0,24}$.  
 By simple facts about root lattice embeddings (see, for example \cite[Section 4.1]{nishi1996}), it turns out that the modular curves are given by the following $ 14$ isomorphic classes of lattices $L = J^\perp/J \in \cG(1,e)$ labeled by its root system:
 \begin{itemize}
     \item  $e=3$: $E_{8}^{2}$, $D_{16}$,
     \item $e=1$:
  $ E_{8}\oplus D_{7} $,
  $E_{7}^{2} \oplus A_1 $, 
  $E_{7}\oplus A_{8}$,
  $D_{15}$,
  $D_{12}\oplus A_{3}$,
  $D_{9}\oplus A_{6}\,(A_8\hookrightarrow A_{15} \subset D_9 \oplus A_{15}) $,
  $A_{15}\,(A_8\hookrightarrow D_9 \subset D_9 \oplus A_{15})$,
  $A_{9}\oplus D_{6}$, 
  $ E_{6}\oplus D_{7}\oplus A_{2}$,
  $A_{8}^{2} $,
  $A_{12}\oplus A_{3}$,
  $ A_{15}\, (A_8 \hookrightarrow A_{24}) $ . 
 \end{itemize}
 Note that there are two isomorphic classes with the same root system $A_{15}$.
Since there is only one orbit of isotropic subgroup of $A_{A_8}\oplus A_{1,1}$ giving primitive embedding of $A_8$ into $H\in \uniL_{0,24}$ and  their direct sum with $A_8$ are contained in the different overlattice, by Lemma \ref{lem:overlattice}, the two classes are different.

(ii). For $n=3$, we have $(A_{\Sigma_3},q_{\Sigma_3}) \cong (A_{E_6 \oplus A_2}, -q_{E_6 \oplus A_2}) $.
 So we consider the primitive embedding of $E_6 \oplus A_2$ and $E_8$ respectively.
 There are actually $ 10$ isomorphic classes of lattices $L = J^\perp/J \in \cG(3,e)$ labeled by its root system:
\begin{itemize}
    \item $e=3$: $E_{8}^{2}$, $D_{16}$.   
    \item $e=1$: 
  $E_{8}\oplus E_6 \oplus A_2$,
  $D_{13}\oplus A_{2} $,
  $D_{10}\oplus A_{5} $,
  $E_{7}\oplus D_{7}$,
  $A_{14}$, 
  $E_{6}^{2} \oplus A_{2}^{2} $,
  $D_{4}\oplus A_{11}$,
  $D_{7}\oplus A_{8}$.
\end{itemize}

(iii).  For $n=2$, we have $(A_{\Sigma_2},q_{\Sigma_2}) \cong (A_{\rL_8}, -q_{\rL_8}) \cong (A_{A_8},q_{A_8}) $ where $\rL_8$ is the even lattice of rank 8 with the following Gram matrix under the basis $\{d_1, \cdots d_7 ,w\}$:
  $$ \left(
\begin{smallmatrix}
-2 &  & 1 &  &  &  &  & 0\\ 
  & -2 & 1  &  &  &  &  & 1\\
 1  & 1 & -2 & 1 &  &  &  & -1\\
     &  & 1 & -2 & 1 &  &  & 0 \\
      &  &  & 1 & -2 & 1 &  & 0\\
       &  &  &  & 1 & -2 & 1 & 0\\
        &  &  &  &  & 1 & -2 &  1\\
 0&1&-1& 0 & 0 & 0 & 1 & -4
\end{smallmatrix}\right) .
  $$
where $\{d_{1} , \cdots,d_{7}\}$ are the bases of $D_{7}$.

By classifying the primitive embedding of $\rL_8$ and $E_8$, we claim there are two isomorphic classes in $\cG(2,3)$ and $9$ classes in $\cG(2,1)$.
The $e=3$ case is as the same as above.
For the $e=1$ case, first note that we have $D_7 = \mathrm{Span}\{d_1,\cdots ,d_7\} \subset \rL_8$ as a primitive sublattice. It's sufficient to consider the embedding of $\rL_8$ into lattices in $\uniL_{0,24}$ with root system $E_8^3,~ E_8\oplus D_{16},~D_{24}, ~D_{12}^{2}, ~ D_{10} \oplus E_7^2, ~ D_9\oplus A_{15},~ D_8^{3}$ and $D_7 \oplus E_6\oplus A_{11}$.
Then the claim follows from the Lemma \ref{n=2 BailyBorel}.

(iv). For type III boundary components, by \cite[Lemma 4.1.2, Proposition 4.1.3]{Scattone1987}, they are bijective to \[
  \text{\{isotropic elements of $A_{\Sigma_n}$\}}/\{\pm 1\}.
\]
For every $n$, there are only $2$ such elements class, one is $0$.
Any isotropic rank $2$ sublattice contains the vector corresponding to $0$, and only those with $J^\perp/ J \in \cG(n,3)$ contain the vector corresponding to the nontrivial discriminant class. 
This completes the proof.

\end{proof}

\subsection{Looijenga's compactifications}

Let $\rD$ be the type IV domain associated to an integral lattice $(L,q)$ of signature $(2,n)$.
Let $\Gamma $ be a congruence arithmetic subgroup of  the stable orthogonal group $\widetilde{O}(L)$.   
Following \cite{Loo03b}, let $\cH$ be a $\Gamma$-invariant hyperplane arrangement, we set $$\rD^{\circ}=\rD-\mathop{\cup} \limits_{H\in \cH} H $$ and define
\begin{itemize}[topsep=5pt, itemsep=5pt, partopsep=5pt]
\item ${\rm PO}(\cH)$: the collection of subspaces $M\subseteq L $ which are intersection of members in $\cH$ meeting $\rD$. 
Denote by $$\pi_M: \PP(L_{\CC})-\PP(M)\longrightarrow \PP(L_\CC/M)$$  the natural projection.  
The projection also defines a natural subdomain $\pi_M\rD^\circ\subseteq \rD^\circ$ ({\it cf}.~\cite[$\S$7]{Loo03b}).
	
\item $\Sigma(\cH)$: the collection of the common intersection of $I^\perp$ and members in $\cH$ containing $I$, where $I$ is a $\QQ$-isotropic line or plane of  $L$.  
\end{itemize}
Then Looijenga's compactification of $\Gamma\backslash \rD^\circ$ can be interpreted as below: let
\begin{equation}\label{looij}
\widehat{\rD}=\rD^{\circ} \cup \coprod\limits_{M\in {\rm PO}(\cH)} \pi_M\rD^\circ \cup \coprod\limits_{V\in \Sigma(\cH)} \pi_{V}\rD^\circ, 
\end{equation}
we define $\overline{\Gamma\backslash \rD}^\cH $ to be the quotient $\Gamma\backslash\widehat{\rD}$, which compactifies $\Gamma\backslash \rD^\circ$ and boundary decomposes into finitely many stratas. 
The birational map $\overline{\Gamma\backslash \rD}^\cH  \dashrightarrow (\Gamma\backslash \rD)^{\ast}$ can be resolved by the following diagram
\begin{equation}\label{BB-Looij}
\xymatrix{
  \widetilde{\Gamma\backslash \rD}^\cH  \ar[d]_{\pi_1} \ar[r]^{\pi_2}\ar[dr]^{\widetilde{\pi}} &    (\Gamma\backslash \rD)^{\Sigma(\cH)}\ar[d]^{\pi_{\cH} }    \\
  \overline{\Gamma\backslash \rD}^\cH   \ar@{.>}[r] &  \Gamma\backslash \rD^{\ast} }
\end{equation}
where $\pi_\cH: (\Gamma\backslash \rD)^{\Sigma(\cH)}  \to  (\Gamma\backslash \rD)^{\ast} $ is the $\QQ$-Cartierization of the hyperplane arrangement  in $\cH$ and  $\pi_i$ are the blow up and blow down respectively. 
Note that when $\cH^{(r)} \neq \emptyset$, i.e. there exists non-empty common intersection of $D$ with $r$ linearly independent hyperplanes in $\cH$, the dimension of the boundary 
\begin{equation} \label{eq:boundary dim}
    \dim (\overline{\Gamma\backslash D}^\cH - \Gamma\backslash D^\circ)\geq  r-1 .
\end{equation}

By \cite[Collary 7,5]{Loo03b}, each component $\pi_V\rD^\circ$ in \eqref{looij} has codimension $\geq 2$, and if $\pi_M \rD^\circ$ has codimension $\geq 2$, the blow up map $\pi_1$ will contract the divisors $\widetilde{\pi}^*(\Gamma \backslash \cup_{H\in \cH} H)$.
In this case, the Looijenga's compactification $ \overline{\Gamma\backslash \rD}^\cH$ can also be interpreted in terms of the algebra of automorphic forms as
\begin{equation}
\begin{split}
       \overline{\Gamma\backslash \rD}^\cH &\cong \proj R(\rD^\circ, \cO(1)|_{\rD^\circ})^\Gamma\\
      &\cong \proj R(\Gamma\backslash \rD^\circ, \lambda|_{\Gamma\backslash \rD^\circ}).
\end{split}
\end{equation}
In particular, take $\cH= \emptyset$, it gives back the Baily--Borel compactification \[
(\Gamma\backslash \rD)^* \cong \proj R(\Gamma\backslash \rD, \lambda).
\]
For $L = \Sigma_1$, $\Gamma = \widetilde{O}(\Sigma_1)$, let $\cH_1 $ be the collection of hyperplanes such that $\Gamma \backslash D^\circ$ is the complement of the Heegner divisors $\bH_u$ and $\bH_h $.
To be specific, we have $\cH_1 = \cH_u \cup \cH_h$ where
\[
    \cH_u = \bigcup_{\substack{v\in\Lambda,\, v^2=0 \\ v\cdot C =1,\, v\cdot E=1 
  }} \pi_{\Sigma_1}(v)^\perp ,\quad\quad 
  \cH_h =    \bigcup_{\substack{v\in\Lambda,\, v^2=0 \\ v\cdot C =2,\, v\cdot E=1
  }} \pi_{\Sigma_1}(v)^\perp 
\]
Then $\overline{\cF}_{T_1}^{\cH_1}$ is the Looijenga's compactification of $\cF_{T_1} - \bH_u \cup \bH_h$.

\begin{lemma}
When $n=1$, the boundary $\overline{\cF}^{\cH_1}_{T_1} - \Gamma \backslash D^\circ $ has codimension $1$.
\end{lemma}

\begin{proof}
Consider the signature $(1,19)$ lattice $N$ given by the span of $C,E, e_1,\cdots, e_{18} $ where 
$e_i\innerp C=2,~ e_i \innerp E=1$ and $e_i\innerp e_j = \delta_{ij}$. 
Since $A_N \cong \ZZ/27\ZZ$, it's easy to show that $N$ can be primitively embedded into the K3 lattice $\Lambda$ by Nikulin 's results (see \cite[Corollary 1.12.3]{Nikulin1979}).
Then $N$ can represent the common intersection of $18$ hyperplanes in $\cH_h$.
This proves the assertion by \eqref{eq:boundary dim}.
\end{proof}

\begin{corollary}  \label{corollary:non-iso}
For $n=1$, the GIT quotient $\overline{\cM_2}^{\rm GIT}$ is not isomorphic to the Looijenga's compactification $\overline{\cF}_{T_1}^{\cH_1}$.
\end{corollary}
\begin{proof}
     Since $\cM_2$ is isomorphic to $\Gamma \backslash D^\circ \cong \cF_{T_1} - \bH_u \cup \bH_h$, this is obtained by comparing the dimension of the boundary of $\overline{\cM}_2^{\rm GIT} - \cM_2$ and $\overline{\cF}_{T_1}^{\cH_1} - \Gamma \backslash D^\circ$.
\end{proof}

\subsection{HKL program for degree \texorpdfstring{$8$}{8} and trielliptic K3 surfaces}
\label{subsec:HKL}
Recall that the Mukai models of degree $8$ K3 surfaces are complete intersections of  three quadrics in $\PP^5$.
Moreover, it was shown in \cite{Greer2015} that Mukai's GIT model $\overline{\cF_8}^{\rm Mukai}\defeq \mathrm{Gr}(3,21)/\!/\SL_6$ is a compactification of the complement of the union of extremal NL-divisors $\rP_{1,1}^8\cup \rP^8_{2,1} \cup \rP^8_{3,1}$.
We first recall that in \cite{Greer}, the Hassett-Keel-Looijenga program for $\cF_8$ is to compare $\overline{\cF_8}^{\rm Mukai}$  and the Baily--Borel compactification $\cF_8^*$ by finding interpolating models in between.
More precisely, consider the log canonical model of $\cF_8^*$ of the form \[
    \overline{\cF}_8(s) = \proj R(\cF_8,\lambda_8 + \Delta_8(s))
\]
where $\lambda_8$ is the Hodge line bundle of $\cF_8$ and $\Delta_8(s)  = \rP^8_{1,1} + \rP^8_{2,1} + s \rP^8_{3,1} $ with a parameter $s\in [0,1]\cap \QQ  $.
Following the vision in  \cite{LO21}, we predict that the birational transformations that occur in \[
\overline{\cF}_8(0) \dashrightarrow \overline{\cF_8}(1)
\]
are almost controlled by the stratification of the support of $\Delta_8(s)$ and $\overline{\cF}_8(s)$ interpolates between $\cF_8^*$ and $\overline{\cF}_8^{\rm{Mukai}}$ in the following way.

Let $\cF^\circ_8$ be the complement of $\rP^8_{1,1}$ and $\rP^8_{2,1}$, then we have
\begin{itemize}
    \item[(0)] the ring of sections \[
    R(\cF_8^\circ ,\lambda_8+s\rP^8_{3,1}|_{\cF_8^\circ}) \cong R(\cF_8, \lambda_8+ \Delta_8(s))
    \]
    is finitely generated for each $s\in [0,1]\cap \QQ$.
    
    \item[(1)] (Semitoric part) $\overline{\cF_8}(0)$ is a Looijenga's compactification associated with the hyperplane arrangements $\cH$ corresponding to the primitive NL divisors $\rP^8_{1,1} $ and $\rP^8_{2,1}$
    and the birational map  \[
    \cF_8^* \dashrightarrow \overline{\cF}_8(0)
    \]
    factors as 
    \begin{equation}
\xymatrix{
  \widetilde{\cF_8}^\cH  \ar[d]_{\pi_1} \ar[r]^{\pi_2}\ar[dr]^{\widetilde{\pi}} &    \cF_8^{\Sigma(\cH)}\ar[d]^{\pi_{\cH} }    \\
 \overline{\cF}_8(0)    &  \cF^*_8 \ar@{-->}[l]
  }
\end{equation}
where $\pi_{\cH}$ is the $\QQ$-Cartierization of $\rP^8_{1,1} \cup \rP^8_{2,1}$ and $\pi_1$ is the divisorial contraction of the $\widetilde{\pi}^*(\rP^8_{1,1} \cup \rP^8_{2,1})$.

\item[(2)] (GIT part) There is an isomorphism \[
    \overline{\cF}^{\rm{Mukai}}_8 \cong \overline{\cF}_8(1)
\]

\item[(3)] the walls of the Mori chamber decomposition of the cone \[
\left\{ (\lambda_8+s\rP^8_{3,1})|_{\cF_8^\circ}) \,\middle| \, s\in \QQ s\geq 0\right\} 
\]
are located at a finite set $\mathbf{Wall}(\cF_8)$ consisting of critical values in $[0,1]\cap \QQ$.

The map \[\overline{\cF}_8(0) \dashrightarrow \overline{\cF}_8(1)\]
admits a factorization by a series of elementary birational transformations (flips and divisorial contractions) between interpolating models $\overline{\cF_8}(s_n, s_{n+1})$ of the form
\[
\begin{tikzcd}
\cF^*_8 \arrow[d, dashed] & {\overline{\cF}_8(0,s_2)} \arrow[ld] \arrow[r, dashed,"\cdots"] & {\overline{\cF}_8(s_{n-1},s_n) } \arrow[rd, "\nu^-_n"] &                       & {\overline{\cF}_8(s_{n},s_{n+1}) } \arrow[ld, "\nu^+_n"'] \arrow[r, dashed,"\cdots"] & {\overline{\cF}_8(s_m,1) } \arrow[d]                   \\
\overline{\cF}_8(0)       &                                                        &                                                        & \overline{\cF}_8(s_n) &                                                                            & \overline{\cF}_8(1)\cong \overline{\cF}_8^{\rm{Mukai}}
\end{tikzcd}
\]
where $s_n\in \mathbf{Wall}$ is the $n$-th critical value.
The center of the birational morphisms are proper transformation of Shimura subvarieties.
\end{itemize}

The study of trielliptic K3 surfaces will provide much help for the above conjecture since there is an natural relation between the HKL program for degree $8$ and type I trielliptic K3 surfaces.
Note that  $\cF_{T_1}$ admits a natural morphism to $\cF_8$
\[
\iota\mcolon \cF_{T_1} =  \Gamma_v \backslash\rD(v^\perp) \rightarrow \Gamma \backslash \rD = \cF_8
\]
as locally symmetric varieties , with image onto the trielliptic divisor $\rP^8_{3,1}$.
Here $\rD$ is the period domain associated with $\Lambda_8 = U^2 \oplus E_8^2 \oplus (-8) $,  $\Gamma = \widetilde{O}(\Lambda_8)$ and $\Gamma_v$ is the stabilizer of the trielliptic vector $v$ such that $\rP^8_{3,1} =\Gamma \backslash \bigcup_{v'\in \Gamma\cdot v} (v')^\perp$.
Following the similar computations as in \cite[Section 5.2, 5.3]{LO19} and more general cases in \cite{Greer}, one have
\[
\iota^*(\rP^8_{1,1}) =0 , \quad  \iota^*(\rP^8_{2,1}) = \bH_u. \quad  \iota^*(\rP^8_{3,1}) = -\lambda_8 +\bH_u + \bH_h 
\]
Then consider the log canonical models of $\cF_{T_1}$ induced by the pullback $\iota^*$
\[
\begin{split}
       \overline{\cF}_{T_1}(s) &\defeq \proj R(\cF_{T_1},(\lambda_8 + \Delta_8(s))|_{\cF_{T_1}}) = \proj R(\cF_{T_1},(1-s)\lambda_{T_1} + (1+s) \bH_u +s \bH_h) \\
    &\cong \proj R(\cF_{T_1}, \lambda_{T_1} + \frac{1+s}{1-s} \bH_u + \frac{s}{1-s} \bH_h)
\end{split}
\]
for $s\in [0, 1)\cap \QQ$.
This will also connect the semitoric compactification and the GIT compactification of $\cF_{T_1}$.
The study of $\overline{\cF}_{T_1}(s)$ will provide the critical walls of $\overline{\cF}_8(s)$ by simple observation.

\begin{remark}
    With the computation and predictions in \cite{Greer}, one can get the critical walls of $\overline{\cF}_{T_1}$.
    For example, the last three are $s = 1,\,\frac{1}{2}$, and $\frac{1}{3}$.
    We will study the HKL program for $\cF_{T_1}$ and $\cF_8$ in the sequel of this paper.
\end{remark}

\section*{Appendix. Classification of some Non-root lattice embeddings}

\subsection*{Notations}
We first fix the following notations. 
\begin{enumerate}
    \item Denote $A_n=\left\langle a_1,a_2,\ldots,a_n\right\rangle$  where $\{a_i\}_{i=1}^n$ is the standard basis.
Let $\alpha_k:=a_k^* $ be the dual basis in $A_n^*$ with 
$\alpha_k^2= \frac{k}{k+1}$.

\item Let $E_n = \langle e_1, \cdots ,e_n\rangle $ with $e_i^2 =-2$, $  e_1\innerp e_4=e_i \innerp e_{i+1}=1$ for $2 \leq i \leq n-1$ and other products equal zero. 
Let $\epsilon_6 = e_6^*$ be the dual basis in $E_6^*$ with $\epsilon_6^2=-\frac{4}{3}$
and $\epsilon_7 = e_7^*$ be the dual basis in $E_7^*$ with  $\epsilon_7^2=-\frac{3}{2}$.

\item Denote $D_n=\left\langle d_1, d_2, \ldots, d_n\right\rangle$ where 
$\{d_i\}_{i=1}^{n}$ is the
standard basis of $D_{n}$
with $d_1\innerp d_3=d_i \innerp d_{i+1}=1$ for $2 \leq i \leq n-1$ and other products equal zero.
Let $\delta_k:=d_k^* $ be the dual basis in $D_n^*$.
 We know $\delta_1, \delta_n$ are the generators for $A_{D_n}$ when $n$ is even, $\delta_1$ is the generator for $A_{D_n}$ when $n$ is odd with $\delta_1^2 \ =-\frac{n}{4}, \, \delta_n^2 =-1, \delta_1 \innerp \delta_n = -\frac{1}{2}$.
 
  \item We denote the corresponding unique unimodular lattice by $M(R)$ for $R$ the root system of some Niemeier lattice in $\uniL_{0,24}$.
 
 \item  The rank $8$ even lattice $\rL_8$ is given by the following Gram matrix under the basis $\{d_1, \cdots d_7 ,w\}$:
  $$ \left(
\begin{smallmatrix}
-2 &  & 1 &  &  &  &  & 0\\ 
  & -2 & 1  &  &  &  &  & 1\\
 1  & 1 & -2 & 1 &  &  &  & -1\\
     &  & 1 & -2 & 1 &  &  & 0 \\
      &  &  & 1 & -2 & 1 &  & 0\\
       &  &  &  & 1 & -2 & 1 & 0\\
        &  &  &  &  & 1 & -2 &  1\\
 0&1&-1& 0 & 0 & 0 & 1 & -4
\end{smallmatrix}\right) 
  $$
  with $(A_{\Sigma_2},q_{\Sigma_2}) \cong (A_{\rL_8}, -q_{\rL_8}) \cong (A_{A_8},q_{A_8}) $.
\end{enumerate}

\begin{lemma}[Eichler's criterion for $D_n$]
\label{lem:Eichler definite}
For any $x \in D_n^*$, $x^2 = \delta_{n}^2 =-1$ if and only if there is $g\in O(D_n)$ such that $g\cdot x =\delta_n$.
\end{lemma}

\begin{proof}
The if part is trivial and we show the only if part below.
Firstly, suppose $(\delta_1-\delta_{n}+\sum_{i=1}^{n}x_id_i )^2=-1$ and we get 
$$
 x_3^2=(x_1-x_3)^2+(1-x_1)^2+(x_3-x_2)^2+
 x_2^2+\sum_{i=3}^{n-1} (x_i-x_{i+1}) +(x_{n}+1)^2.
$$
From \[
\begin{split}
    (x_1-x_3)^2+(1-x_1)^2 &\geq \frac{1}{2}(x_3-1)^2  \\
  (x_3-x_2)^2+x_2^2 &\geq \frac{1}{2} x_3^2  \\
  \sum_{i=3}^{n-1} (x_i-x_{i+1}) +(x_{n}+1)^2 &\geq | x_3+1 |
\end{split}
  \]
we get a contradiction.
Since $\delta_1 -\delta_n$ and $\delta_1$ are connected by the reflection with respect to $\delta_n$, one see $x$ must be of the form $\delta_{n}+\sum_{i=1}^{n}x_i d_i$.

Secondly, suppose $( \delta_{n}+\sum_{i=1}^{n}x_i d_i )^2=-1$.
We get 
$$
 x_3^2=(x_1-x_3)^2+x_1^2+(x_3-x_2)^2+
 x_2^2+\sum_{i=3}^{n} (x_i-x_{i+1}) +(x_{n}+1)^2.
$$
  From $$(x_1-x_3)^2+x_3^2 \geq \frac{1}{2}{x_3^2}$$
  $$(x_3-x_2)^2+x_2^2 \geq \frac{1}{2}{x_3^2}$$
   $$\sum_{i=3}^{n} (x_i-x_{i+1})^2 +(x_{n}+1)^2 \geq |x_3+1 |$$
we know $x_3 \in \left\{-2,-1,0\right\} $.
 One can easily check that any two vectors of the form $\delta_{n}+\sum_{i=1}^{n}x_id_i$  with  $x_3 \in \left\{-2,-1,0\right\} $
  can be related by reflections 
  of $d_i$.

\end{proof}

\begin{remark}\label{rmk:Eichler ADE}
    A weaker result also holds for $(A_{11},\alpha_{11})$, $(2\alpha_{15}+A_{15}, 2a_{15}+\alpha_{15})$, $(E_6, \epsilon_6)$ and $(E_7,\epsilon_7)$ by applying the similar methods as in the proof of the above Lemma.
    More precisely, for $(\Lambda, v\in \Lambda^*)$ in the above list, if 
    \[
    (v+u)^2 = v^2 \text{~and~} u\in \Lambda,
    \]
    we have $v+u = g\cdot v$ for some $g\in O(\Lambda)$.
\end{remark}

\begin{lemma}  \label{n=2 BailyBorel}
Considering the primitive embedding of $\rL_8$ in
a unimodular lattice of rank $24$, there are $9$ isomorphic classes of the orthogonal complement uniquely determined by their root system:  
 \begin{itemize}
   
    \item $E_7 \oplus E_8 $ \quad   $(\rL_8 \hookrightarrow M(E_8^{\oplus 3}))$
     
    \item $A_1 \oplus D_{14} $\quad $(\rL_8 \hookrightarrow M(E_8 \oplus D_{16}))$
 
    \item $E_{8} \oplus A_8 $  \quad $(\rL_8 \hookrightarrow M(D_{16}\oplus E_8  ))$

    \item $A_4 \oplus D_{11}$ \quad  $(\rL_8 \hookrightarrow M(D_{12}^{\oplus 2}))$

    \item $E_7 \oplus E_6 \oplus A_2$ \quad $(\rL_8 \hookrightarrow M( D_{10}\oplus E_7^{\oplus 2}))$

    \item $A_1 \oplus A_1 \oplus A_{13}$ \quad $(\rL_8 \hookrightarrow M(D_9 \oplus A_{15} ))$   

    \item $D_7 \oplus D_7$ \quad $(\rL_8 \hookrightarrow M(D_8 ^{\oplus 3}))$     

     \item $D_8 \oplus A_6$ \quad $(\rL_8 \hookrightarrow M(D_8 ^{\oplus 3}) )$      

    \item $D_5 \oplus A_{10}$  \quad$(\rL_8 \hookrightarrow M( D_{7}\oplus E_6 \oplus A_{11}))$    
    
\end{itemize}

 \end{lemma}

 \begin{proof}
As noted in Theorem \ref{thm:Baily--Borel}, $\rL_8$ can only be primitively embedded into $M(R)$ for \[
R\in\{E_8^3,~ E_8\oplus D_{16},~D_{24}, ~D_{12}^{2}, ~ D_{10} \oplus E_7^2, ~ D_9\oplus A_{15},~ D_8^{3},~ D_7 \oplus E_6\oplus A_{11} \}.
\]
(1). First we deal with the cases where the root system $R$ contains a direct summand of $D_n$ type with $n \neq 16 $ in a relatively unified way.  
In these cases, one can recover  $M(R)$ as a sublattice of  $ D^*_{n}\oplus T$.
For a primitive embedding $\tau \mcolon \rL_8 \hookrightarrow M(R)$, we may assume $\tau(d_i) = d_i \in D_{n}$ for $ 1\leq i\leq 7$ up to an automorphism of $M(R)$.
Let $p$ be the projection of $D_n^*\oplus T$ onto $D_n^*$.
Write
\[p ( \tau (w)) =\sum _{i=1}^{n} x_{i} d_{i} +k\delta_1+l \delta_n \]
where $x_{i} \in \mathbb{Z}$ and $k,l \in \mathbb{Z}/2\mathbb{Z}$ since $2\delta_i =0 \pmod {D_n}$ for possible $R$.

By solving the Diophantine equations given by the Gram matrix of $\rL_8$, one can assume  $k=-1$.
When $n\geq 8$, we get 
\begin{align*}
     p(\tau (w))&=-\delta_{1} +l \delta_n+c(d_{1} +d_{2}) + (1+2c) \sum_{i=3}^7 d_i +( 2+2c) d_{8} + \sum_{j=9}^n 
 x_j d_j.
\end{align*}
We denote the vector on the right hand-side by $w_{(x_9,\cdots,x_n)}$ for $c$ and $l$ fixed.
Write
$ p(\tau (w))^2 =g_{1} +g_{2}$ where 
\[
\begin{split}
    g_{1} &=\delta_1^2-2x_{1} -2\sum_{i=1}^7 x_i^2 -x_{8}^{2} +2( x_{1} x_{3} +\sum_{i=2}^7 x_{i} x_{i+1} ) = -2+\delta_1^2 -2c,\\
  g_2 &= -\sum_{i=8}^{n-1}( x_i -x_{i+1})^{2} -x_{n}^{2}+2x_nl. 
\end{split}
\]
Besides, we have $g_{2} \leq 0$ when $l=0$.
In below, we calculate the orthogonal complement of $\rL_8\hookrightarrow M(R)$ case by case.
\begin{itemize}

\item $M(D_{24})$:

This is the \textbf{only case} that $\rL_8$ admits \textbf{no} primitive embedding.
Suppose on the contrary we have primitive embedding $\tau \colon  \rL_8 \hookrightarrow M(D_{24})$. 
Note that $$
M(D_{24}) \cong \ZZ\cdot \delta_1 +D_{24} \subset (D_{24})^*.
$$  
Hence $l=0$.
 We have 
 \[
 \begin{split}
     -4  &= p( \tau(w) )^2 = -2c-8 - \sum_{i=8}^{23}(x_i - x_{i+1})^2-x_{24}^2 \\
     &\leq -2c-8 - |x_8| = -2c-8 -|2+2c| \\
     &< -4
 \end{split}
 \] 
which is a contradiction.

\item  $M(D_{12} \oplus D_{12})$:

In this case, we show $\rL_8^\perp$ is unique up to isomorphism for any $\rL_8 \xrightarrow{\tau} M(D_{12}^2)$.   
Write 
\[
    M(D_{12} \oplus D_{12}) \cong D^{2}_{12} + \ZZ\cdot(\delta_1,\delta_{12})  + \ZZ\cdot (\delta_{12},\delta_1).
\]
Let $D_7 \subset \rL_8$ be embedded into the first factor $\D_{12}$.
Then we can write $\tau(w)$ as \[
\tau(w) = (p(\tau(w), l\delta_1 -\delta_{12} +v ),\quad v\in D_{12}.
\]
If $l=1$, then one get \begin{equation}\label{eq:D12}
     -4 \leq p(\tau(w))^2 \leq g_1+1 -|x_8 - 1|= -2c-4- |2c+1| =-3.
\end{equation}
  
Since the second coordinate of $\tau(w)$ can't have zero norm, one have $p(\tau(w))^2 =-3$ and $(\delta_1-\delta_{12}+v)^2 = \delta_{12}^2=-1$.
This contradicts the Lemma \ref{lem:Eichler definite} above.
Therefore, one can assume $l=0$ and the second coordinate is $\delta_{12}$ by the Lemma \ref{lem:Eichler definite} again.

Therefore for fixed $c$, the embedding $\tau$ is uniquely determined by the orbit of $w_{(x_9,\cdots,x_{12})}$ under $\Gamma_{D_7} \subset O(M(D_{12}^2))$ that fix $D_7$.
Moreover, we have $p(\tau (w))^2=-3$ and  $c \in \left\{-1,-2,-3 \right\} $ by the simple estimate as in \eqref{eq:D12}.

 Firstly, for each fixed $c$, the corresponding orthogonal complement ${\rL_8^{\perp}}(c)$ is independent of 
 the choice of $(x_8,\cdots,x_{12})$.
 This is because one can relate any $w_{(x_8,\cdots,x_{12})}$ by a composition of reflections $\rho_{d_i}$ associated with the roots $d_9,d_{10},d_{11},d_{12}$ which fix $D_7$.
 For instance, considering $\rL_8^\perp(-2)$, one sees the embeddings are connected as
 \[
 \begin{split}
      w_{(-2,-2,-2,-2,-1)}\xrightarrow{\rho_{d_{11}}}   w_{(-2,-2,-2,-1,-1)} \\
      \xrightarrow{\rho_{d_{10}} }
w_{(-2,-2,-1,-1,-1)} \xrightarrow{\rho_{d_9} }
w_{(-2,-1,-1,-1,-1)} .
 \end{split}
 \]

Secondly, we show $\rL_8^\perp(-1)\cong \rL_8^\perp(-2) \cong \rL_8^\perp(-3) $. 
Note that the orthogonal complement of $D_7$ in $D_{12}$ is $\left\langle d_8',d_9,\ldots, d_{12} \right\rangle$
where $d_8'=d_1+d_2+2(\sum_{i=3}^8 d_i)+d_9$.
For $c=-1$, $p(\tau (w))=-\delta_1-\sum_{i=1} ^7 d_i $ and we denote the vector by $v_1$.
 For $c=-3$, $p(\tau (w))=-\delta_1-\sum_{i=1} ^7 d_i -4d_8 -3d_9-2d_{10}-d_{11}$ and we denote the vector by $v_3$.
 Notice that $d_8' \innerp v_3=1=-d_8'\innerp v_1$.
 Then one can construct a isomorphism $\psi$ which induces 
 a isomorphism between $\rL_8^\perp(-1)$  and $\rL_8^\perp(-3)$:
\[
 \begin{split}
     \psi \mcolon M(D_{12} \oplus D_{12} )  &\longrightarrow M( D_{12} \oplus D_{12}) \\
     ( x,z) &\mapsto (-x,z)
 \end{split}
\]
 The isomorphsim $\psi$ also gives the isomorphsim $\rL_8^\perp(-1) \cong \rL_8^\perp(-2)$ for the same reason.
 And one easily sees that the root lattice of the unique orthogonal complement is  $R(\rL_8^\perp(-1))=A_4 \oplus D_{10}$.

\item  $M(D_9 \oplus A_{15}) $:

 In this case,  $l=0$ as
$$
M(D_9 \oplus A_{15}) \cong \ZZ\cdot (\delta_1,2\alpha_{15}+a_{15})+D_9 \oplus A_{15} \subset (D_9 \oplus A_{15})^*.
$$  
Given embedding $\tau$, by estimating $p(\tau(w))$ as in the above case, we get  
$c \in \left\{-1,-2 \right\} $ as 
$-\frac{17}{4}-2c=g_{1} \geq -\frac{9}{4} $ and $((2+2c)-x_9)^2+x_9^2 = -(2+2c)$. 
The same method in the case $D_{12} \oplus D_{12}$ 
implies $\rL_8^\perp(-1) \cong \rL_8^\perp(-2)$.
The embedding is given by mapping $w$ to $(-\delta_1-\sum_{i=1} ^7 d_i,2\alpha_{15}+a_{15})$ as the vector $(2\alpha_{15}+a_{15})$ is unique in the sense of Remark \ref{rmk:Eichler ADE}.
The root system $R(\rL_8^\perp) \cong A_1 \oplus A_1 \oplus A_{13}$.


\item $ M(D_{10}\oplus E_7^{\oplus 2})$:

Write\[
M(D_{10}\oplus E_7^{\oplus 2}) \cong \ZZ\cdot(\delta_1,0,0) + \ZZ\cdot(0,\epsilon_7,0) + D_{10}\oplus E_7^{\oplus 2}.
\]
 In this case, $l=0$ and we have the unique embedding isomorphism class.
 By the same argument, we get $-3 \leq c \leq -1$.
The embedding can be given by mapping $w$ to $(-\delta_1-\sum_{i=1} ^7 d_i,-\epsilon_7,0)$.
The root system of the orthogonal complement  is $A_2  \oplus E_6 \oplus E_7 $.

\item $ M(D_{8}^{\oplus 3})$

  In this case,
there are two different kinds of embeddings corresponding to  two possible values of $l$.
 We recover  $M(D_{8} \oplus D_{8} \oplus D_8)$, as the sublattice $ D^*_{8}\oplus D^*_{8} \oplus D^*_{8}$ spanned by $D^{\oplus 3}_{8}$,
     $(\delta_1,\delta_{8},\delta_{8})$,
     $(\delta_8,\delta_{1},\delta_{8}) $ and $(\delta_8,\delta_{8},\delta_{1})$.
     
     If $l=0$, one can see that $p (\tau (w))^2 =-4c^2-10c-8$. 
    Then $c=-1$ and $w$ is mapped to $(-\sum_{i=1}^{7}d_{i}-\delta_1, \delta_8,\delta_8)$
    with $R(\rL_8^\perp) \cong D_7 \oplus D_7$. 

     If $l=1$, by calculation we know 
    $ p( \tau (w)) =-\sum_{i=1}^{7}d_{i} -\delta_1 +\delta_8$,
 In this case, $w$ is mapped to $(-\sum_{i=1}^{7}d_{i} -\delta_1 +\delta_8, \delta_1 -\delta_8,0)$
    with $R(\rL_8^\perp) \cong D_8 \oplus A_6$.

\item $M(D_{7}\oplus E_6 \oplus A_{11})$:
Consider a primitive embedding $\tau\mcolon \rL_8 \hookrightarrow M(D_{7}\oplus E_6 \oplus A_{11})$.
Write $$
M(D_{7}\oplus E_6 \oplus A_{11}) \cong \ZZ\cdot(\delta_1,\epsilon_6,\alpha_{11})+D_{7}\oplus E_6 \oplus A_{11}\subset (D_{7}\oplus E_6 \oplus A_{11})^*.
$$ 
Suppose $\tau(D_7) = D_7$ the identity map and $p(\tau(w)) = k\delta_1 + \sum_{i=1}^7 x_i  d_i$. 
Assuming $k =0$ or $-1$ and solving the equations given by the Gram matrix, we get the unique solution $x_i=k=-1$.
Then we get the unique isomorphism classes of primitive embedding by Remark \ref{rmk:Eichler ADE} with $R(\rL_8^\perp) \cong D_5 \oplus A_{10}$.

\end{itemize}

(2). Now we deal with the remaining two cases $E_8^{\oplus 3}$ and $M(E_8\oplus D_{16}) =E_8 \oplus D_{16}^+ $.     
For the case of $E_8^{3}$, there is only one primitive sublattice $\tau \mcolon \rL_8 \hookrightarrow E_8^2 \subset E_8^3$.
Note that  $E_8\cong D_8 + \ZZ \cdot \delta_1  = \mathrm{Span}\{ \delta_1,d_1,\cdots, d_7\}$, with $ e_2 \leftrightarrow \delta_1,\, e_3 \leftrightarrow d_1, e_1 \leftrightarrow d_2,  e_{i+1} \leftrightarrow d_i$ for $3\leq i \leq 7$,
One can assume that the primitive embedding sends $D_7\hookrightarrow D_8 \subset E_8$  into the first $E_8$ factor with the above identification. 
Since $w^2 =-4$, $\tau(w) $ must lie in the first $E_8$ plus another one, say the second $E_8$ up to a permutation of the three factors.
Then the unique embedding can be given by mapping $w$ to $ (-\delta_1- \sum_{i=1}^7 d_i,v,0)$ where $v$ is a $(-2)$ vector in the second copy of $E_8$.
In this case $R(\rL_8^\perp) \cong E_8\oplus (A_1\hookrightarrow E_8)^\perp = E_8 \oplus E_7$.
Note that the genus of $A_8$ contains two classes.
The primitive embedding of $\rL_8 \hookrightarrow E_8^2$ actually gives the another one which contains $E_7$.

For the case of $E_8 \oplus D_{16}^+ $, there are two different primitive sublattices that are isomorphic to $\rL_8$. 
Firstly, a similar calculation as in $M(D_{12}^2)$ shows that $\rL_8$ can be primitively embedded into $D_{16}^+ \cong D_{16} + \ZZ \cdot \delta_1$, hence also $E_8\oplus D^+_{16}$.
Note that the orthogonal complement is in the genus of $A_8$ which contains two classes. The other one contains the root sublattice $E_7$ hence it cannot lie in $D^+_{16}$. 
This shows that the orthogonal complement of $\rL_8$ in $D^+_{16}$ is unique, hence also the complement in $E_8\oplus D^+_{16}$.
The corresponding root system is $E_8\oplus A_8$.
Another kind of embedding can be given as follows.
Notice that $E_8\cong D_8 + \ZZ\cdot \delta_1^{(8)}$, here we specify $\delta_1^{(8)}\in D_8^*$.
Then the unique embedding isomorphism class $\rL_8 \hookrightarrow E_8 \oplus D_{16}^+$ is given by letting $w \mapsto (-\delta_1^{(8)}- \sum_{i=1}^7 d_i,v)$ where $v\in D_{16}^+$ is a $(-2)$-vector. 
The corresponding root system is $(A_1\rightarrow D_{16})^\perp \cong A_1 \oplus D_{14}$.

 \end{proof}

\bibliographystyle{alpha}
\bibliography{ref.bib}
\end{document}